\documentclass[12pt,a4paper]{article}
\usepackage{amssymb,amsmath,amsfonts,latexsym} 
\usepackage{hyperref}
   
\usepackage{mathtools}
\usepackage{multirow,array}
\usepackage{subcaption}
\usepackage{graphicx}
\usepackage{epstopdf}
\usepackage{float} 
\usepackage{color, xcolor}

\usepackage{amsfonts}
\usepackage{amsmath}
\usepackage{amsbsy}
\usepackage{theorem}
\usepackage{graphicx}

\newtheorem{lemma}{Lemma}

\begin{document}

\sf
\title{\bfseries Manifold functional multiple regression model with LRD error term
}
\author{Diana  P. Ovalle--Mu\~{n}oz and M. Dolores Ruiz--Medina}
\date{}

 \maketitle
\begin{abstract}
This paper considers the problem of manifold functional multiple  regression with functional response, time--varying scalar regressors,  and functional error term displaying  Long Range Dependence (LRD) in time. Specifically, the error term is given  by a manifold   multifractionally integrated functional
time series (see, e.g.,  Ovalle--Mu\~noz \&  Ruiz--Medina \cite{Ovalle23}).
  The manifold is defined by a connected and compact two--point homogeneous space. The functional regression parameters have support in the manifold.   The Generalized Least--Squares (GLS) estimator of the vector  functional regression parameter is computed, and its asymptotic properties are analyzed under a totally specified and misspecified model scenario. A   multiscale residual correlation analysis  in the  simulation study  undertaken illustrates the empirical distributional properties of the errors at different spherical resolution levels.
\end{abstract}


\section{Introduction}
\label{sec:1}

There exists an extensive literature on functional linear regression given  its wide interesting applications in several scientific fields. Among many other contributions, we refer to environmental applications such as air pollution studies (see, e.g., Acal et.al. \cite{acal:2022}; \'Alvarez--Li\'ebana \& Ruiz--Medina \cite{ALVRuiz19}; Olaya--Ochoa, Ovalle--Mu\~noz \& Urbano--Leon \cite{Olaya20};  and references therein),
biomechanics applications involving  modeling of human movements (see, e.g., in Acal \& Aguilera \cite{acal:2023}; Helwig et. al. \cite{Helwig:2016}), and epidemiological and brain applications (see, e.g., Aristizabal, Giraldo \& Mateu \cite{Aristizabal:2019};
Yao, M\"uller \& Wang \cite{yao:2005}). The main directions of contributions in the above cited references  focused on  predictive analysis of the response  based on functional  regressors in both cases, when response is scalar or functional,  as well as  in the case where the regression parameters are functions in the  parametric framework.  In particular, dimension reduction techniques based on functional principal component analysis, constitute a major topic in this literature. Most of the regularization techniques proposed are based on projection into suitable finite--dimensional subspaces, whose dimension depends on the functional sample size. Also suitable weighting operators are considered for embedding the unbounded inverse autocovariance operator of the regressors into the space of bounded linear operators (see, e.g., Bosq \cite{Bosq2000};   Cardot, Mas \& Sarda \cite{Cardot07}; Crambes \&  Ma \cite{Crambes13},   Mas \cite{Ma99}, among others). The main subject of these papers is to remove  the ill--posed nature of the associated inverse estimation problem due to unbounded inversion of the autocovariance operator of the functional covariates. Indeed, during the last twenty years several contributors have worked on several numerical proposals for approximation of the slope function, and its asymptotic analysis. The corresponding residual correlation analysis has also been addressed under different scenarios, mainly including the cases of independent or weak--dependent data (see, e.g., Chiou \& M\"uller \cite{Chiou07}, Shen \& Xu \cite{Shen}, among others).

The case of  independent functional observations was extensively developed during  the first decade of 2000. In particular, regression analysis and inference from a sample of independent and identically distributed functional random variables have been considered in  the papers by  Crambes, Kneip \&  Sarda  \cite{Crambes09};
Cuevas,   Febrero \&   Fraiman  \cite{Cuevas02}; Febrero--Bande,  Galeano \&  Gonzalez--Manteiga \cite{FebreroBande15}  (see  also Cuevas  \cite{Cuevas14} for an overview). The   kernel  formulation of the  regression parameters is usually  adopted in the literature of parametric linear  regression with  functional response and regressors (see, for example,   Chiou,  M\"uller \&  Wang \cite{Chiou04}; Ruiz--Medina (2011) \cite{RuizMedina11}; Ruiz--Medina  (2012a) \cite{RuizMedina12a}; Ruiz--Medina (2012b) \cite{RuizMedina12b},  and  references therein).
An extensive overview, with further references on functional regression approaches, including the case of functional response and regressors, can be as well  found in Morris  \cite{Morris15}. See also the   monograph by  Hsing \& Eubank \cite{HsingEubank15}, where several  functional analytical tools are introduced, for estimating random elements in function spaces.

An extended formulation of the previous results to the case of weak--dependent functional observations is provided in the context of functional time series, and in general, of weak--dependent processes. That is the case of approaches based on
 the concept of $L^{r}$--$m$--approximability to modeling  the temporal dependence in the regression functional errors (see, for example,  Horv\'ath \&  Kokoszka \cite{HorvathandKokoszka}).  A central topic in this book is the analysis of functional data,  displaying  dependent  structures    in time and space. Also, in the framework of weakly dependent functional time series  models,  supporting inference on stochastic processes, in particular, in a state space  framework, several functional regression approaches have been adopted for functional  prediction (see, e.g., \'Alvarez--Li\'ebana et.al. \cite{Alvarez22};  Guillas \cite{Guillas2002};  H\"ormann \& Kokoszka \cite{Kokoszka10}; Horv\'ath \& Kokoszka \cite{HorvathandKokoszka}; Kara--Terki \& Mourid \cite{Kara2016}; 	Kokoszka \& Reimherr \cite{Kokoszka2013}).  Confidence bands, kernel  and parametric functional time series estimation lead to the analysis of outstanding  problems like adaptive bandwidth selection, dimension reduction, change--point analysis,  and functional principal component estimation, among others (see, e.g.,  Berkes, Horv\'ath, \& Rice \cite{Berkes16}; Dette, Kokot \& Aue \cite{Dette20}; H\"ormann \& Kidzinski \cite{HOKID15};  H\"ormann \& Kokoszka \cite{Kokoszka13};  Horv\'ath, Rice \& Whipple \cite{Horvath2016};  Zhang et.al. \cite{Zhang11}).

A fixed effect approach in Hilbert spaces is adopted in Ruiz--Medina \cite{RuizMedina16}, for FANOVA analysis under dependent errors. For simple regression, with explanatory variable taking  values in some abstract space of functions, the rate of convergence of the mean squared error of the functional version of the Nadaraya--Watson kernel estimator is derived in  Benhenni,  Hedli--Griche \& Rachdi \cite{Benhenni},  when the errors are represented by a  stationary short or long memory process. An alternative approach, based on Autoregressive Hilbertian processes of order 1 (ARH(1)) error tem, is presented in Ruiz--Medina, Miranda \& Espejo \cite{RuizMD18} for multiple regression, and in \'Alvarez--Li\'ebana \& Ruiz--Medina \cite{Alvarez17} for fixed effect models including the case of circular domains.

Classical texts like   the book by Ramsay  \&  Silverman  \cite{Ramsay05} have introduced the basic analytical and statistical tools for inference on stochastic processes based on Functional Data Analysis (FDA). The book by  Ferraty \&  Vieu  \cite{Ferraty06} also constitutes a benchmark  in the literature on  nonparametric functional statistics.   In this framework, one can mention, among others,  the contributions by Aneiros--P\'erez \cite{Aneiros06}, Aneiros--P\'erez \& Vieu \cite{Aneiros08}, and Ferraty et al.  \cite{FerratyGoia13}, applying, in particular,  the Projection Pursuit Regression principle    in  the approximation of the regression function, for the case of a functional predictor and a scalar response   (see also Ferraty \& Vieu (2006) \cite{Ferraty06}; Ferraty \& Vieu (2018) \cite{Ferraty11}).  In the nonparametric setting, Ferraty,  Keilegom \& Vieu \cite{FerratyKeilegom12}  derive  a kernel type estimator of the regression operator, and its  asymptotic normality is proved, for the case of functional response and predictor.
Goia \&  Vieu \cite{GoiaVieu15} adopt a semiparametric approach,  in  a two--terms Partitioned Functional Single Index Model. Several papers on  topics related to  statistical analysis of high--dimensional data, including functional regression, from the parametric, semiparametric  and nonparametric FDA frameworks,  can be found in the Special Issue  by Goia \& Vieu \cite{GoiaVieu16}.

Recently,  an attempt to extend spectral analysis of  functional time series   to the context of  LRD functional sequences has been presented  in Ruiz--Medina \cite{RuizMedina2022}, covering, in particular,  some examples of the LRD funtional time series family analyzed by Li, Robinson \& Shang \cite{LiRobinsonShang19} in the temporal domain. In particular, Li, Robinson \& Shang \cite{LiRobinsonShang19}  applies Functional Principal Component Analysis (FPCA)  based on the long--run covariance function,  for the consistent estimation of the dimension and the orthonormal functions  spanning the dominant subspace, where the projected curve process displays the largest dependence range. Fractionally integrated functional  autoregressive moving averages processes  constitute an  interesting example (see Li, Robinson \& Shang \cite{LiRobinsonShang19}).
The  multifractional version of this process family  can be analyzed under the modeling framework introduced in Ruiz--Medina \cite{RuizMedina2022}.
Indeed, in that paper, one can see the  important  advantages that
the application   of harmonic analysis entails  in this more general context.  Particularly, under stationary in time, the temporal dependence range can be approximated from the behavior in a neighborhood of zero frequency of the spectral density operator family  at different spatial resolution levels. Moreover, a more flexible modeling framework can be introduced in this setting, allowing the representation of  long, intermediate or short range dependence, according to the interval where the pure point  or continuous spectrum of the   long--memory operator lies (see, e.g., Ovalle--Mu\~noz \& Ruiz--Medina \cite{Ovalle23}).

In the case of  $\mathbb{M}_{d}$  being a  connected  and compact two--point homogeneous space   with topological dimension $d,$ the  invariance of a kernel   with respect to the group of isometries of $\mathbb{M}_{d}$ allows its spectral  diagonalization in terms of a fixed  orthogonal basis, the eigenfunctions of the Laplace Beltrami operator on  the Hilbert space $H=L^{2}\left(\mathbb{M}_{d},d\nu\right).$
   Here, $d\nu $ denotes the measure  induced by the probabilistic invariant measure on  the connected component of the group of isometries of $\mathbb{M}_{d}.$  In the functional time series framework, one can find recent contributions on modeling, estimation and asymptotic analysis of weak--dependent Hilbert--valued processes in the sphere, which constitutes a well--known example of connected and compact two--point homogeneous space  (see, e.g., Caponera \& Marinucci, \cite{CaponeraMarinucci}; Caponera \cite{Caponera21}).   Formulation of alternative limit results for sojourn measures of LRD  spherical--cross time random fields can be found in  Marinucci, Rossi \& Vidotto \cite{MarinucciRV}.  The reader can see the starting results for manifold analysis based on connected and compact two--point homogeneous spaces  in the preliminary  contributions on this subject obtained in Ma \& Malyarenko \cite{MaMalyarenko}, for second--order mean--square continuous elliptically contoured random fields  on  $\mathbb{M}_{d}.$ Several motivating applications of this modeling framework are given, for example,  in  Alegr\'{i}a et al. \cite{Alegria21}, on  bayesian multivariate spherical random field modeling  from a  bivariate spatial data set from two 2019 NCEP/NCAR Flux reanalyses, as well as, in  Marinucci  \&  Peccati \cite {MarinucciPeccati11}, and Leonenko, Nanayakkara \& Olenko  \cite{Leonenko21}, in  relation to  Cosmic Microwave Background   (CMB) evolution modeling and data analysis.

The present paper  applies the results in Ruiz--Medina \cite{RuizMedina2022} and Ovalle--Mu\~noz \& Ruiz--Medina \cite{Ovalle23}  to  introduce a new infinite--dimensional linear  model  family in the context of multiple functional regression models  evaluated in a Hilbert space, with finite--dimensional design matrix, and functional response and regression parameters. This paper focuses in the special case of functional error term being an  LRD functional time series in the family introduced in Ruiz--Medina \cite{RuizMedina2022}.  Thus,   compact Riemannian manifold--based residual correlation analysis is addressed here from  functional regression parameter estimation when the error term displays LRD. This analysis is achieved from the previous theoretical and  numerical results obtained in Ovalle--Mu\~noz \& Ruiz--Medina \cite{Ovalle23}  for LRD functional time series in connected  and compact two--point homogeneous spaces. Thus, this paper derives  an extended formulation of the linear  functional models introduced in Ruiz--Medina \cite{RuizMedina16}; and \'Alvarez--Li\'ebana \& Ruiz--Medina \cite{Alvarez17}, beyond the
  weak--dependent and  Euclidean settings.

The  outline of the paper is as follows. Section \ref{secprel}  presents some preliminary elements on the spectral analysis of LRD manifold functional time series. Section \ref{secFANOVA} introduces our multiple functional regression setting in a parametric framework. The  generalized least--squares  estimator of the manifold   functional regression  parameter vector is then computed.  In Section \ref{secsimulation}, a simulation study is undertaken to illustrate the finite--sample and asymptotic properties of the theoretical and empirical functional response predictor.  In particular, the effect of the spectral properties  of the LRD operator, characterizing the temporal dependence range of the  functional error term, on the precision and variability of these response predictors is analyzed.  Some final comments and open research lines are discussed in Section \ref{secconclusions}.

\section{Preliminaries}

\label{secprel}
In the regression residual correlation analysis  achieved in this paper in the spectral domain, under  LRD manifold functional time series error term,
 invariance of the involved kernels   with respect to the group of isometries of $\mathbb{M}_{d}$ plays a crucial role.  Indeed, these kernels are diagonalized by the eigenfunctions $\{ S_{n,j}^{d},\ j=1,\dots,\delta(n,d), \ n\in \mathbb{N}_{0}\}$ associated  with the  eigenvalues  $\{ \lambda_{n}=-n\varepsilon(n\varepsilon+\alpha+\beta+1), \ n\in \mathbb{N}_{0}\}$ of  the Laplace Beltrami operator  $\Delta_{d}$ on $L^{2}(\mathbb{M}_{d},d\nu ,\mathbb{R})$ (see, e.g.,   Cartan \cite{Cartan1927} and Ma  \&   Malyarenko   \cite{MaMalyarenko}, for more details on Lie Algebra based  approach).

 We first formulate the addition formula applied in the context of connected and compact two--point homogeneous spaces.
    \begin{lemma}
\label{cor:adformula}
 (See \cite[Theorem 3.2.]{Gine1975} and \cite[p 455]{Andrews99})
For every $n\in \mathbb{N}_{0},$ the following addition formula holds:
\begin{equation}
\sum_{j=1}^{\delta(n,d)}S_{n,j}^{d}(\mathbf{x})S_{n,j}^{d}(\mathbf{y})=\frac{\delta(n,d)}{\omega_{d}}R_{n}^{\alpha,\beta }\left(\cos(d_{\mathbb{M}_{d}}(\mathbf{x},\mathbf{y}))\right),\quad \mathbf{x},\mathbf{y}\in \mathbb{M}_{d}.
\label{af}
\end{equation}
\noindent Here, $\omega_{d}=\int_{\mathbb{M}_{d}}d\nu(\mathbf{x}),$
and  $\delta(n,d)$ denotes the dimension of the eigenspace $\mathcal{H}_{n}$ associated with the  eigenvalue  $\lambda_{n}=-n\varepsilon(n\varepsilon+\alpha+\beta+1)$ of
the Laplace Beltrami operator, which   is given, for every $n\in \mathbb{N}_{0},$  by
\begin{equation}
\label{ecdeltand}
\delta(n,d)=\frac{(2n+\alpha+\beta +1)\Gamma (\beta +1)\Gamma (n+\alpha +\beta +1)\Gamma (n+\alpha +1)}{\Gamma (\alpha +1)\Gamma (\alpha +\beta +2)\Gamma (n+1)\Gamma (n+\beta +1)}.
\end{equation}
\noindent Furthermore,
 $R_{n}^{\alpha,\beta}\left(\cos(d_{\mathbb{M}_{d}}(\mathbf{x},\mathbf{y}))\right)
=\frac{P_{n}^{\alpha,\beta}\left(\cos(d_{\mathbb{M}_{d}}(\mathbf{x},\mathbf{y}))\right)}{P_{n}^{\alpha,\beta}\left(1\right)},$ with $P_{n}^{\alpha,\beta}$ denoting the Jacobi polynomial of degree $n\in\mathbb{N}_{0},$ with parameters $\alpha$ and $\beta,$ involved in the definition of $\delta(n,d)$ (see equation (\ref{ecdeltand})).
\end{lemma}

  Let   $X=\{ X(\mathbf{x},t),\ \mathbf{x}\in \mathbb{M}_{d},\ t\in \mathbb{T}\}$  be  a zero--mean, stationary in time, and isotropic  in space mean--square continuous Gaussian, or elliptically contoured, spatiotemporal random field  on the basic probability space $(\Omega,\mathcal{A},P),$  with covariance function
$C(d_{\mathbb{M}_{d}}(\mathbf{x},\mathbf{y}),t-s)= E\left[X(\mathbf{x},t) X(\mathbf{y},s)\right],$ for $\mathbf{x},\mathbf{y}\in \mathbb{M}_{d},$ and $t,s\in \mathbb{T}.$ Here, $\mathbb{T}$ denotes the temporal domain, which can be $\mathbb{Z}$ or $\mathbb{R}.$ Under the conditions of   Theorem  4   in  Ma \& Malyarenko  \cite{MaMalyarenko}, the covariance function $C(d_{\mathbb{M}_{d}}(\mathbf{x},\mathbf{y}),t-s)$ admits the following diagonal series expansion:
\begin{eqnarray}&&C(d_{\mathbb{M}_{d}}(\mathbf{x},\mathbf{y}),t-s)=
  \sum_{n\in \mathbb{N}_{0}}
B_{n}(t-s)\sum_{j=1}^{\delta (n,d)}S_{n,j}^{d}(\mathbf{x})S_{n,j}^{d}(\mathbf{y})\nonumber\\
&&=
\sum_{n\in \mathbb{N}_{0}} \frac{\delta (n,d)}{\omega_{d}}  B_{n}(t-s)R_{n}^{(\alpha, \beta )}\left(\cos\left(d_{\mathbb{M}_{d}}(\mathbf{x},\mathbf{y})\right)\right),\ \mathbf{x},\mathbf{y}\in \mathbb{M}_{d}, \ t,s \in \mathbb{T}.\nonumber\\ \label{klexpc2}
\end{eqnarray}

Let now  $X=\{ X(\mathbf{x},t),\ \mathbf{x}\in \mathbb{M}_{d},\ t\in [0,T]\}$  be the restriction to the interval $[0,T]$ of a zero--mean, stationary in time, and isotropic in space,  mean--square continuous Gaussian, or elliptically contoured, spatiotemporal random field  on the basic probability space $(\Omega,\mathcal{A},P),$  with covariance function $C(d_{\mathbb{M}_{d}}(\mathbf{x},\mathbf{y}),t-s)$ admitting the diagonal  expansion  (\ref{klexpc2}). The following lemma provides the orthogonal expansion of $X=\{ X(\mathbf{x},t),\ \mathbf{x}\in \mathbb{M}_{d},\ t\in [0,T]\}$ in terms of the eigenfunctions of the Laplace Beltrami operator (see Theorem 1 in the Supplementary Material in Ovalle--Muñoz \& Ruiz--Medina \cite{Ovalle23}).
\begin{lemma}
\label{lemmat5ma}
Let $X=\{ X(\mathbf{x},t),\ \mathbf{x}\in \mathbb{M}_{d},\ t\in [0,T]\}$  be the restriction to the interval $[0,T]$ of a zero--mean, stationary in time,  and isotropic  in space, mean--square continuous Gaussian, or elliptically contoured, spatiotemporal random field  on the basic probability space $(\Omega,\mathcal{A},P),$  with covariance function (\ref{klexpc2}) satisfying  the conditions in  Theorem 4 in Ma \& Malyarenko \cite{MaMalyarenko}, with
\begin{equation}\sum_{n\in \mathbb{N}_{0}}B_{n}(0)\delta (n,d)<\infty.\label{eqtrace}
\end{equation}
\noindent  Then,  the following orthogonal expansion holds for   random field
 $X$
 \begin{eqnarray}
&&X(\mathbf{x},t)\underset{\mathcal{L}^{2}_{\widetilde{H}}(\Omega,\mathcal{A},P)}{=}\sum_{n\in \mathbb{N}_{0}} \sum_{j=1}^{\delta (n,d)}V_{n,j}(t)S_{n,j}^{d}(\mathbf{x}),\quad \mathbf{x}\in \mathbb{M}_{d},\ t\in [0,T],
\label{klexp}
\end{eqnarray}

\noindent where $\mathcal{L}^{2}_{\widetilde{H}}(\Omega,\mathcal{A},P)=L^{2}(\Omega\times\mathbb{M}_{d}\times[0,T],P(d\omega )\otimes d\nu\otimes dt),$ with
  $\widetilde{H}=L^{2}(\mathbb{M}_{d}\times[0,T],
 d\nu\otimes dt).$ Here,
 $\{V_{n,j}(t),\ t\in  [0,T], \ j=1,\dots, \delta(n,d), \  n\in \mathbb{N}_{0}\}$ is a sequence of centered   random processes on $[0,T]$ given by
 \begin{equation}
 V_{n,j}(t)=\int_{\mathbb{M}_{d}}X(\mathbf{y},t)S_{n,j}^{d}(\mathbf{y})d\nu(\mathbf{y}), \  j=1,\dots, \delta(n,d),\ n\in \mathbb{N}_{0},
 \label{fcA}
 \end{equation}
 \noindent in the mean--square sense.
\end{lemma}

Assume  that $\mathbb{T}=\mathbb{Z},$ and that the map $$\widetilde{X}_{t}:(\Omega ,\mathcal{A})\longrightarrow \left(L^{2}(\mathbb{M}_{d},d\nu ,\mathbb{R}),\mathcal{B}(L^{2}(\mathbb{M}_{d},d\nu ,\mathbb{R}))\right)$$ \noindent  is   measurable, with  $\widetilde{X}_{t}(\mathbf{x}):=X(\mathbf{x},t)$ for every $t\in \mathbb{T}$ and $\mathbf{x}\in \mathbb{M}_{d}.$  Here,  $\mathcal{B}(L^{2}(\mathbb{M}_{d},d\nu ,\mathbb{R})) $ denotes the Borel  $\sigma$--algebra on $L^{2}(\mathbb{M}_{d},d\nu, \mathbb{R})$
(i.e., the smallest $\sigma$--algebra containing  the collection of all open subsets of  $L^{2}(\mathbb{M}_{d},d\nu, \mathbb{R})$).
By previous assumptions on $X,$  $\left\{\widetilde{X}_{t}, \ t\in \mathbb{Z}\right\}$   then defines a manifold stationary functional time series. In particular, $E\left[\widetilde{X}_{t}\right]=0,$ and   $\sigma_{\widetilde{X}}^{2}=E\left[\|X_{t}\|^{2}_{L^{2}(\mathbb{M}_{d},d\nu,\mathbb{R})}\right]=
E\left[\|\widetilde{X}_{0}\|^{2}_{L^{2}(\mathbb{M}_{d},d\nu,\mathbb{R})}\right]=\|R_{0}\|_{L^{1}(H)},$ for every $t\in \mathbb{Z}.$ By $L^{1}(L^{2}(\mathbb{M}_{d},d\nu,\mathbb{R}))$ we denote  the space of trace or nuclear operators on $L^{2}(\mathbb{M}_{d},d\nu,\mathbb{R}).$
The second--order structure of $\left\{\widetilde{X}_{t}, \ t\in \mathbb{Z}\right\}$ is  characterized by the family of covariance operators  $\left\{\mathcal{R}_{t}, \ t\in \mathbb{Z}\right\}$  given by, for all $h,g\in L^{2}(\mathbb{M}_{d},d\nu,\mathbb{R}),$
\begin{eqnarray}
&&
\mathcal{R}_{t}(g)(h)= E[\widetilde{X}_{s+t}(h)\widetilde{X}_{s}(g)]= E\left[\left\langle
\widetilde{X}_{s+t},h\right\rangle_{L^{2}(\mathbb{M}_{d},d\nu,\mathbb{R})}\left\langle \widetilde{X}_{s},g\right\rangle_{L^{2}(\mathbb{M}_{d},d\nu,\mathbb{R})}\right]\nonumber\\
\label{covstatfunct2}
\end{eqnarray}
\noindent with respective kernels
$$r_{t}(\mathbf{x},\mathbf{y})=E[\widetilde{X}_{s+t}\otimes \widetilde{X}_{s}](\mathbf{x},\mathbf{y}),\quad
\forall  \mathbf{x},\mathbf{y} \in \mathbb{M}_{d},\   t,s\in \mathbb{Z}.$$

We now introduce some preliminary elements about spectral analysis of  functional time series,  based on the spectral density operator family, and the periodogram operator, computed from the functional discrete Fourier transform (fDFT) (see, e.g.,  Panaretos \& Tavakoli \cite{Panaretos13}, for the weak--dependent case; Ruiz--Medina \cite{RuizMedina2022}, Ovalle--Mu\~{n}oz \& Ruiz--Medina \cite{Ovalle23}, for the strong--dependent case).
 The functional Fourier transforms of the elements of the family of covariance operators introduced in (\ref{covstatfunct2}) are here defined  in the Hilbert--Schmidt operator norm (see, e.g.,  Ovalle--Mu\~noz \& Ruiz--Medina \cite{Ovalle23}). That is, the family of spectral density operators  $\left\{\mathcal{F}_{\omega},\ \omega \in [-\pi, \pi]\right\}$  characterizing the second--order structure of the functional time series $\{\widetilde{X}_{t},\ t\in \mathbb{Z}\}$  in the spectral domain is given by,
for each $\omega \in [-\pi,\pi]\backslash \{0\}:$
  \begin{equation}
  \mathcal{F}_{\omega}
  \underset{\mathcal{S}(L^{2}(\mathbb{M}_{d},d\nu,\mathbb{C}))}{=}
  \frac{1}{2\pi} \sum_{t\in \mathbb{Z}}\exp\left(-i\omega t\right)
\mathcal{R}_{t},\label{sdo2}
\end{equation}
\noindent where $\underset{\mathcal{S}(L^{2}(\mathbb{M}_{d},d\nu,\mathbb{R}))}{=}$ denotes the identity in the norm of the Hilbert--Schmidt operators. In particular, this identity does not require the summability of the series of trace norms of the elements of the covariance operator family. Equivalently,  short--memory is not assumed in our functional time series  modeling framework (see, e.g., Panaretos \& Tavakoli \cite{Panaretos13}, Ruiz--Medina  \cite{RuizMedina2022} and Ovalle--Mu\~noz \& Ruiz--Medina \cite{Ovalle23}).

The fDFT $\widetilde{X}^{(T)}_{\omega }(\cdot)$ of the  curve data  is defined as
\begin{equation}\widetilde{X}^{(T)}_{\omega }(\cdot)
\underset{L^{2}(\mathbb{M}_{d},d\nu,\mathbb{C})}{=} \frac{1}{\sqrt{2\pi
T}}\sum_{t=1}^{T}\widetilde{X}_{t}(\cdot )\exp\left(-i\omega t\right), \quad
\omega\in [-\pi ,\pi],\label{fDFT}\end{equation}\noindent
 where $\underset{L^{2}(\mathbb{M}_{d},d\nu,\mathbb{C})}{=}$ denotes the equality
in $L^{2}(\mathbb{M}_{d},d\nu,\mathbb{C})$ norm, with $L^{2}(\mathbb{M}_{d},d\nu,\mathbb{C})$ being the complex version of the Hilbert space  $L^{2}(\mathbb{M}_{d},d\nu,\mathbb{R}).$  Note that $\widetilde{X}^{(T)}_{\omega }(\cdot)$ is a random  element in the space   $L^{2}(\mathbb{M}_{d},d\nu,\mathbb{C}),$
since $$E\left[\|\widetilde{X}_{\omega}^{(T)}\|_{L^{2}(\mathbb{M}_{d},d\nu,\mathbb{C})}\right]\leq \frac{1}{\sqrt{2\pi
T}}\sum_{t=1}^{T}E\|\widetilde{X}_{t}(\cdot )\|_{\widetilde{H}}<\infty.$$

For $\omega\in [-\pi ,\pi],$ the periodogram operator  $p_{\omega }^{(T)}=\widetilde{X}_{\omega}^{(T)}\otimes
\overline{\widetilde{X}_{\omega}^{(T)}} $ is defined  from the fDFT. Its mean is given by convolution of the  F\'ejer kernel $F_{T}(\omega )=\frac{1}{T}\sum_{t=1}^{T}
\sum_{s=1}^{T}\exp\left(-i(t-s)\omega \right)$ with the spectral density operator.
That is,
  \begin{eqnarray}E[p_{\omega }^{(T)}]&=&E[\widetilde{X}_{\omega}^{(T)}
 \otimes \widetilde{X}_{-\omega}^{(T)}]=  \frac{1}{2\pi }\sum_{u=-(T-1)}^{T-1}\exp\left(-i\omega u\right)
\frac{(T-|u|)}{T}\mathcal{R}_{u}\nonumber \\&=&\int_{-\pi}^{\pi} F_{T}(\omega - \xi)
\mathcal{F}_{\xi} d\xi,\quad T\geq 2.\nonumber \end{eqnarray}

\section{Multiple  Functional regression model}
\label{secFANOVA}

This section introduces our manifold multiple functional regression model in a parametric framework, under an LRD  functional time series modelling of the error term. Specifically, let us consider the following functional observation model in time:  For $t=1,\dots,N,$ and $j=1,\dots,p,$
\begin{equation}
\label{ecmodelo}
Y_{t}(\mathbf{x})=\sum_{j=1}^{p}X_{t,j}\beta_{j}(\mathbf{x})+\varepsilon_{t}(\mathbf{x}),\quad \mathbf{x}\in \mathbb{M}_{d},
\end{equation}
\noindent where  $\beta_{j}\in L^{2}(\mathbb{M}_{d},d\nu ,\mathbb{R}),$ and   $X_{t,j}\in \mathbb{R},$ for  $j=1,\dots,p,$ and $t\in \mathbb{Z}.$ Here,
$\{\varepsilon_{t},\ t\in \mathbb{Z}\}$ defines a stationary  LRD  functional time series with values in the space  $L^{2}(\mathbb{M}_{d},d\nu,\mathbb{R})$  as
 introduced in the previous section. Note that
 $\mathbf{X}=(X_{t,j})_{(t,j)\in \{1,\dots, N\}\times \{1,\dots, p\}}$ defines the $N\times p$ design matrix     (see also Ruiz--Medina \cite{RuizMedina16}). Equation (\ref{ecmodelo}) can be equivalently expressed in vectorial form as follows:
 \begin{equation}
\label{ecmodelov}
\mathbf{Y}(\mathbf{x})=\mathbf{X}\boldsymbol{\beta}(\mathbf{x})+\boldsymbol{\varepsilon}(\mathbf{x})\quad \mathbf{x}\in \mathbb{M}_{d},
\end{equation}
\noindent where $\mathbf{Y}(\mathbf{x})=[Y_{1}(\mathbf{x}),Y_{2}(\mathbf{x}),\dots,Y_{N}(\mathbf{x})]^{T},$  $\boldsymbol{\beta}(\mathbf{x})=[\beta_{1}(\mathbf{x}),\dots, \beta_{p}(\mathbf{x}) ]^{T},$
 and $\boldsymbol{\varepsilon}(\mathbf{x})=[\varepsilon_{1}(\mathbf{x}),\varepsilon_{2}(\mathbf{x}),\dots,\varepsilon_{N}(\mathbf{x})]^{T},$
 for every $\mathbf{x}\in \mathbb{M}_{d}.$

From equation (\ref{ecmodelov}), the second--order structure of the error term $\boldsymbol{\varepsilon}=\{\varepsilon_{t},\ t\in \mathbb{Z}\}$ admits an  infinite--dimensional matrix representation as follows:

\[
\begin{split}
\mathbf{R}_{\boldsymbol{\varepsilon}\boldsymbol\varepsilon}&=E\left[\boldsymbol{\varepsilon}(\cdot)\boldsymbol{\varepsilon}^{T}(\cdot)\right]\\
 &=\begin{bmatrix}
E\left[\varepsilon_{1}(\cdot)\otimes\varepsilon_{1}(\cdot)\right] & E\left[\varepsilon_{1}(\cdot)\otimes\varepsilon_{2}(\cdot)\right] & \cdots & E\left[\varepsilon_{t_1}(\cdot)\otimes\varepsilon_{N}(\cdot)\right]\\
E\left[\varepsilon_{2}(\cdot)\otimes\varepsilon_{1}(\cdot)\right]& E\left[\varepsilon_{2}(\cdot)\otimes\varepsilon_{2}(\cdot)\right]& \cdots & E\left[\varepsilon_{2}(\cdot)\otimes\varepsilon_{N}(\cdot)\right]\\
\vdots & \vdots & \vdots \\
E\left[\varepsilon_{N}(\cdot)\otimes\varepsilon_{1}(\cdot)\right] & E\left[\varepsilon_{N}(\cdot)\otimes\varepsilon_{2}(\cdot)\right] & \cdots & E\left[\varepsilon_{N}(\cdot)\otimes\varepsilon_{N}(\cdot)\right]
\end{bmatrix}\\
&=\begin{bmatrix}
R_{0} & R_{1} & \cdots & R_{N-1}\\
R_{1} & R_{0}& \cdots & R_{N-2}\\
\vdots & \vdots & \vdots \\
R_{N-1}& R_{N-2}& \cdots & R_{0}.
\end{bmatrix},
\end{split}
\]
\noindent in terms of the elements $\left\{ R_{0}, R_{1},\dots, R_{N-1}\right\}$ of the covariance operator family of the functional time series $\left\{\widetilde{X}_{t}, \ t\in \mathbb{Z}\right\},$ whose values at $t=1,\dots,N,$ are here denoted as $\varepsilon_{1},\dots,\varepsilon_{N}.$
Note that  the functional entries of $\mathbf{R}_{\boldsymbol{\varepsilon}\boldsymbol\varepsilon}$ admit the diagonal series expansion introduced in equation (\ref{klexpc2}) under the conditions assumed in Theorem  4   in  Ma \& Malyarenko  \cite{MaMalyarenko}.
In the subsequent development we will consider the orthogonal representation of the functional regression parameters $\beta_{j},$ $j=1,2,\dots,p,$ and of  the observed functional values of the response variable $Y_{t},$ $t=1,2,\dots,N,$  with respect to  the orthonormal basis  
 $\left\{S_{n,k}^{d}(\mathbf{x}),\ k=1,\dots, \delta (n,d),\ n\in \mathbb{N}_{0}\right\}$ of eigenfunctions of the Laplace Beltrami operator. That is, for every $j=1,\dots,p,$
\begin{equation}
\label{ecseriesbeta}
\beta_{j}(\mathbf{x})=\sum_{n\in\mathbb{N}_{0}}\beta_{n,j}\sum_{k=1}^{\delta(n,d)}S_{n,k}^{d}(\mathbf{x}),\quad \forall \mathbf{x}\in \mathbb{M}_{d}.
\end{equation}
\noindent Furthermore, under the conditions   in  Lemma \ref{lemmat5ma}, let us consider the orthogonal expansion of the observed values of the response in the basis  \linebreak   $\left\{S_{n,k}^{d}(\mathbf{x}),\ k=1,\dots, \delta (n,d),\ n\in \mathbb{N}_{0}\right\},$ for every $\mathbf{x}\in 
\mathbb{M}_{d},$ and for each $t\in \mathbb{Z},$
\begin{equation}
\label{ecexpansionres}
\begin{split}
Y_{t}(\mathbf{x})&=\sum_{j=1}^{p}
\sum_{n\in\mathbb{N}_{0}}\sum_{k=1}^{\delta(n,d)}Y_{n,k}(t)S_{n,k}^{d}(\mathbf{x})=\sum_{j=1}^{p}\sum_{n\in\mathbb{N}_{0}}
\sum_{k=1}^{\delta(n,d)}\left[X_{t,j}\beta_{n,j}+V_{n,k}(t)\right]S_{n,k}^{d}(\mathbf{x}),
\end{split}
\end{equation}
\noindent  where convergence holds pointwise on $\mathbb{M}_{d}$  in the mean--square sense for each fixed $t=1,\dots,N,$ and in the space
 $\mathcal{L}^{2}_{\widetilde{H}}(\Omega,\mathcal{A},P)=L^{2}(\Omega\times\mathbb{M}_{d}\times[0,N],P(d\omega )\otimes d\nu\otimes dt),$ under  conditions   of Lemma \ref{lemmat5ma}.

\subsection{GLS functional  parameter estimation}

The $[L^{2}(\mathbb{M}_{d},d\nu, \mathbb{R})]^{p}-$valued GLS functional parameter estimator of  \linebreak $\boldsymbol{\beta}=\left[\beta_{1},\beta_{2},\dots,\beta_{p}\right]^{T}$ is computed from projection into  the orthonormal basis
 $\{ S_{n,k}^{d},\ k=1,\dots,\delta(n,d),\ n\in \mathbb{N}_{0}\}$ of eigenfunctions of the Laplace Beltrami operator
$\Delta_{d}$ on $L^{2}(\mathbb{M}_{d},d\nu,\mathbb{R}),$  in the spirit of the approach adopted in Ruiz--Medina \cite{RuizMedina16}.
Specifically,  for each $n\in \mathbb{N}_{0},$ one can define the matrix

$$
\boldsymbol{\Lambda}_{n}=\begin{bmatrix}
B_{n}(0) & \cdots & B_{n}(N-1)\\
\vdots & \vdots & \vdots\\
B_{n}(N-1) & \cdots & B_{n}(0)
\end{bmatrix},
$$
\noindent  and the respective vector response and error term projections
$$\mathbf{Y}_{n}=\left[Y_{n}(1),Y_{n}(2),\dots,Y_{n}(N)\right]^{T} \   \mbox{and} \   \boldsymbol{\varepsilon}_{n}=\left[V_{\mathbf{n}}(1),V_{\mathbf{n}}(2),\dots,V_{\mathbf{n}}(N)\right]^{T}$$
\noindent   with $\mathbf{n}=(n_{1},\dots n_{\delta(n,d)}),$ into the eigenspace $\mathcal{H}_{n}$ of the Laplace Beltrami operator,    as well as the corresponding  vector of projections $\boldsymbol{\beta }_{n}=\left[\beta_{n,1},\dots,\beta_{n,p}\right]^{T}$ of the functional regression parameter.

\noindent Hence, the minimizer $\widehat{\boldsymbol{\beta}}$   of the mean quadratic loss function $L,$ given by
\begin{eqnarray}
L&=& \left\|\mathbf{Y}-\mathbf{X}\boldsymbol{\beta }\right\|_{\mathbf{R}^{-1}_{\boldsymbol{\varepsilon}\boldsymbol\varepsilon}}
 =\sum_{n\in \mathbb{N}_{0}}\left[\mathbf{Y}_{n}-\mathbf{X}\boldsymbol{\beta}_{n}\right]^{T}\boldsymbol{\Lambda}_{n}^{-1}
\left[\mathbf{Y}_{n}-\mathbf{X}\boldsymbol{\beta}_{n}\right]\nonumber\\
&=&\sum_{n\in \mathbb{N}_{0}}\left\|\mathbf{\varepsilon}_{n}\right\|^{2}_{\boldsymbol{\Lambda}_{n}^{-1}},
\end{eqnarray}
\noindent has Fourier coefficients $\widehat{\boldsymbol{\beta}}_{n}$ defined as
\begin{equation}
\label{ecbetamcg}
\widehat{\boldsymbol{\beta}}_{n}=\left[\hat{\beta}_{n,1},\hat{\beta}_{n,2},\dots,\hat{\beta}_{n,p}\right]^{T}
=(\mathbf{X}^{T}\boldsymbol{\Lambda}^{-1}_{n}\mathbf{X})^{-1}\mathbf{X}^{T}\boldsymbol{\Lambda}_{n}^{-1}\mathbf{Y}_{n},\quad n\in\mathbb{N}_{0}.
\end{equation}

Thus, the theoretical predictor  is defined by
\begin{equation}
\widehat{\mathbf{Y}}=\mathbf{X}\widehat{\boldsymbol{\beta}}
\label{thpred}
\end{equation}

\subsection{Distributional characteristics  of GLS estimator}
\label{secpropertbetaest}
 For every $n\in \mathbb{N}_{0},$ $\hat{\boldsymbol{\beta}}_{n}$ is an unbiased estimator of $\boldsymbol{\beta}_{n}$ since

\begin{equation}
\begin{split}
E\left[ \hat{\boldsymbol{\beta}}_{n}\right]=& E\left[(\mathbf{X}^{T}\boldsymbol{\Lambda}_{n}^{-1}\mathbf{X})^{-1}\mathbf{X}^{T}\boldsymbol{\Lambda}_{n}^{-1}\mathbf{Y}_{n}\right]\\
=& (\mathbf{X}^{T}\boldsymbol{\Lambda}_{n}^{-1}\mathbf{X})^{-1}\mathbf{X}^{T}\boldsymbol{\Lambda}_{n}^{-1}E\left[\mathbf{Y}_{n}\right]\\
=&(\mathbf{X}^{T}\boldsymbol{\Lambda}_{n}^{-1}\mathbf{X})^{-1}\mathbf{X}^{T}\boldsymbol{\Lambda}_{n}^{-1}\mathbf{X}\boldsymbol{\beta}_{n}\\
=&\boldsymbol{\beta}_{n}.
\end{split}
\label{ub}
\end{equation}

Hence, from (\ref{ub}),  for every $\mathbf{x}\in \mathbb{M}_{d},$

\[
\begin{split}
&E\left[\hat{\boldsymbol{\beta}}(\mathbf{x})\right]=E\left[\left(\sum_{n=0}^{\infty}\hat{\beta}_{n,1}\sum_{k=1}^{\delta(n,d)}S^{d}_{n,k}(\mathbf{x}),\dots,\sum_{n=0}^{\infty}\hat{\beta}_{n,p}\sum_{k=1}^{\delta(n,d)}S^{d}_{n,k}(\mathbf{x})\right)^{T}\right]\\
&=\left(\sum_{n=0}^{\infty}E\left[\hat{\beta}_{n,1}\right] \sum_{k=1}^{\delta(n,d)}S^{d}_{n,k}(\mathbf{x}),\dots,\sum_{n=0}^{\infty}E\left[\hat{\beta}_{n,p}\right] \sum_{k=1}^{\delta(n,d)}S^{d}_{n,k}(\mathbf{x})\right)^{T}
\\
&=\left(\sum_{n=0}^{\infty}\beta_{n,1}\sum_{k=1}^{\delta(n,d)}S^{d}_{n,k}(\mathbf{x}),\dots,\sum_{n=0}^{\infty}\beta_{n,p}\sum_{k=1}^{\delta(n,d)}S^{d}_{n,k}(\mathbf{x})\right)^{T}\\
&= \boldsymbol{\beta}(\mathbf{x}).
\end{split}
\]

Finally, as it is well known, since for every $n\in \mathbb{N}_{0},$  \begin{equation}\hat{\boldsymbol{\beta}}_{n}=\boldsymbol{\beta }_{n}+(\mathbf{X}^{T}\boldsymbol{\Lambda}_{n}^{-1}\mathbf{X})^{-1}\mathbf{X}^{T}\boldsymbol{\Lambda}_{n}^{-1}\boldsymbol{\varepsilon}_{n},
\label{equb}
\end{equation}
\noindent we have

\begin{eqnarray}
\mbox{Var}\left[\hat{\boldsymbol{\beta}}_{n}\right]&=&
E\left[ \left(\hat{\boldsymbol{\beta}}_{n}-\boldsymbol{\beta}_{n}\right)^{T}\left(\hat{\boldsymbol{\beta}}_{n}-\boldsymbol{\beta}_{n}\right)\right]
\nonumber\\
&=&
(\mathbf{X}^{T}\boldsymbol{\Lambda}_{n}^{-1}\mathbf{X})^{-1}\mathbf{X}^{T}\boldsymbol{\Lambda}_{n}^{-1}\boldsymbol{\Lambda}_{n}
\boldsymbol{\Lambda}_{n}^{-1}\mathbf{X}(\mathbf{X}^{T}\boldsymbol{\Lambda}_{n}^{-1}\mathbf{X})^{-1}
\nonumber\\
&=&(\mathbf{X}^{T}\boldsymbol{\Lambda}_{n}^{-1}\mathbf{X})^{-1}.
\label{vgls}
\end{eqnarray}

\subsection{Functional spectral based  plug--in estimation of the functional regression parameter}
\label{secmiss}
This section presents a plug--in GLS estimation methodology when the second order structure of the error term is unknown. In our case,  the entries of the matrix sequence $\left\{\boldsymbol{\Lambda}_{n},\ n\in \mathbb{N}_{0}\right\}$ are misspecified. The approach presented is based on   the estimation of such entries
in the spectral domain by minimum contrast,  under a semiparametric modelling framework, applying the methodology introduced in  Ovalle--Mu\~noz \& Ruiz--Medina \cite{Ovalle23}, and Ruiz--Medina \cite{RuizMedina2022}. Specifically, we will assume that the function elements of the Fourier transform sequence
$$f_{n,\theta}(\omega)=\int_{-\pi}^{\pi}\exp(-i\omega t)\widehat{B}_{n,\theta}(t)dt, \  \omega \in [-\pi,\pi],\ n\in \mathbb{N}_{0},$$
\noindent admit the following semiparametric modelling: For every $n\in \mathbb{N}_{0},$
\begin{eqnarray}f_{n,\theta}(\omega)&=&
B_{n}^{\eta}(0)M_{n}(\omega )\left[4(\sin(\omega /2))^{2}\right]^{-\alpha (n,\theta)/2},\ \theta \in \Theta,\ \omega \in [-\pi,\pi],\nonumber\\
\label{eqsmc1}\end{eqnarray}
\noindent where $\alpha(n,\theta),$ $M_{n}(\omega),$ and $B_{n}^{\eta}(0)$ are the eigenvalues of the LRD  operator
  $\mathcal{A}_{\theta },$ of the Short Range Dependence  (SRD) spectral family $\left\{\mathcal{M}_{\omega},\ \omega \in [-\pi,\pi]\right\},$ and of the autocovariance operator $R^{\eta}_{0}$ of the manifold white noise innovation process $\eta $ involved in the definition of the error term $\varepsilon,$ respectively.
Hence,  the elements of the spectral density operator family $\left\{\mathcal{F}_{\omega,\theta },\ \omega \in [-\pi,\pi]\right\}$
introduced in equation (\ref {sdo2}) are respectively  approximated by  $$\left\{\mathcal{F}_{\omega,\widehat{\theta }_{N} },\ \omega \in [-\pi,\pi]\right\},$$ \noindent  with $\widehat{\theta }_{N}$ denoting the minimum contrast estimator of parameter $\theta $ based on a functional sample of size $N$  (see equations (3.8)--(3.16) in Ovalle--Mu\~noz \& Ruiz--Medina \cite{Ovalle23}, or equations (5.1)--(5.19) in Ruiz--Medina \cite{RuizMedina2022}). Fourier transform inversion formula
$$\widehat{B}_{n,\hat{\theta}_{N}}(t)=\int_{-\pi}^{\pi}\exp(i\omega t)f_{n,\widehat{\theta}_{N}}(\omega)d\omega,\quad n\in \mathbb{N}_{0},$$
\noindent  is then applied to obtain the estimated  matrix sequence

\[
\widehat{\boldsymbol{\Lambda}}_{n,\widehat{\theta}_{N}}=\left(\begin{bmatrix}
\widehat{B}_{n,\widehat{\theta}_{N}}(0) & \cdots & \widehat{B}_{n,\widehat{\theta}_{N}}(N-1)\\
\vdots & \vdots & \vdots\\
\widehat{B}_{n,\widehat{\theta}_{N}}(N-1) & \cdots & \widehat{B}_{n,\widehat{\theta}_{N}}(0)
\end{bmatrix}\right),\quad n\in \mathbb{N}_{0}.
\]
\noindent The plug--in GLS estimators  of the Fourier coefficients of the functional multiple regression vector  parameter $\boldsymbol{\beta },$ with respect to the orthonormal basis of eigenfunctions of the Laplace Beltrami operator,  admit the following expression:
\[
\widehat{\boldsymbol{\beta}}_{n,\widehat{\theta}_{N}}=\left(\mathbf{X}^{T}
\widehat{\boldsymbol{\Lambda}}^{-1}_{n,\widehat{\theta}_{N}}\mathbf{X}\right)^{-1}\mathbf{X}^{T}
\widehat{\boldsymbol{\Lambda}}^{-1}_{n,\widehat{\theta}_{N}}\mathbf{Y}_{n},\quad n\in \mathbb{N}_{0},
\]
\noindent and the corresponding  plug--in predictor is then computed as
\begin{equation}
\widehat{\mathbf{Y}}_{n,\widehat{\theta}_{N}}=\mathbf{X}\widehat{\boldsymbol{\beta}}_{n,\widehat{\theta}_{N}}.
\label{emppred}
\end{equation}

\section{Asymptotic properties of the GLS functional parameter estimator}

Model (\ref{ecmodelo}) can be interpreted as  a particular case of the  multiple  functional   regression model introduced in equation (1) in
Ruiz--Medina, Miranda \& Espejo \cite{RuizMD18}, when the Hilbert space $H=L^{2}\left(\mathbb{M}_{d},d\nu, \mathbb{R}\right)$ is considered.  Specifically,
the functional regression model (\ref{ecmodelo}) can be isometrically identified with multiple functional regression model in equation (1) in
Ruiz--Medina, Miranda \& Espejo \cite{RuizMD18}, in the
   special case of diagonal kernel regressors with one non--null coefficient or eigenvalue.

Theorem 1 in Ruiz--Medina, Miranda \& Espejo \cite{RuizMD18} provides the conditions for the asymptotic normality of the GLS $\widehat{\boldsymbol{\beta }}$ estimator of the  vector functional regression parameter $\boldsymbol{\beta },$ in the case where the second--order structure of the functional error term $\varepsilon$ is known. In our case, when the function sequence $\left\{f_{n}(\omega ), \ \omega \in [-\pi,\pi],\ n\in \mathbb{N}_{0}\right\}$ is totally specified, and characterizes the pure  point spectra of the elements of the spectral density operator family \linebreak $\left\{ \mathcal{F}_{\omega },\ \omega \in [-\pi,\pi]\right\}$ of the error term.  Furthermore, in this case, Assumptions A1--A6 in Ruiz--Medina, Miranda \& Espejo \cite{RuizMD18} ensure the strong--consistency of   the GLS estimator $\widehat{\boldsymbol{\beta }}$  of $\boldsymbol{\beta },$ as given in their    Theorem 2.
Finally,   under Assumptions $\widetilde{A5}-\widetilde{A6}$   in Ruiz--Medina, Miranda \& Espejo \cite{RuizMD18}, in the case of unknown $\left\{f_{n}(\omega ), \ \omega \in [-\pi,\pi],\ n\in \mathbb{N}_{0}\right\}$ sequence, Proposition 1 in  Ruiz--Medina, Miranda \& Espejo, \cite{RuizMD18},  adapted to our time--varying finite--dimensional design   matrix $\mathbf{X},$ provides the strong--consistency, in the norm of the space  $[L^{2}\left(\mathbb{M}_{d},d\nu, \mathbb{R}\right)]^{p},$ of  the ordinary least--squares estimator of the  functional regression parameter $\boldsymbol{\beta }.$  From Theorem 2 in Ruiz--Medina \cite{RuizMedina2022} (see also Theorem 2 in Ovalle--Mu\~noz \& Ruiz--Medina \cite{Ovalle23}), applying Fourier transform inversion formula, the   consistency in the sense of the integrated weighted mean square error,  measured in the  norm of the space of Hilbert--Schmidt operators $\mathcal{S}(L^{2}\left(\mathbb{M}_{d},d\nu, \mathbb{C}\right))$ on $L^{2}\left(\mathbb{M}_{d},d\nu, \mathbb{C}\right),$   of the plug--in GLS estimator of $\boldsymbol{\beta }$  follows.
Hence, the weak--consistency of the plug--in GLS estimator of $\boldsymbol{\beta }$ in the norm of  $\mathcal{S}(L^{2}\left(\mathbb{M}_{d},d\nu, \mathbb{C}\right))$ also holds.
\section{Simulation study}
\label{secsimulation}
The simulation study undertaken in this section illustrates the  sample properties of the theoretical and plug--in functional response predictors introduced in equations (\ref{thpred}) and (\ref{emppred}), respectively. The numerical results are displayed under  misspecified, and totally specified model of the functional error term.
Specifically, the  projected  (into the basis of eigenfunctions of the Laplace Beltrami operator) empirical   mean--quadratic and absolute errors are computed.  Their empirical distributional characteristics are  visualized in terms of their histograms.  Thus, at each eigenspace of the Laplace Beltrami operator,  the projected residuals are analyzed, under decreasing and increasing sequence of the eigenvalues of the LRD operator $\mathcal{A}_{\theta }$ in equation
(\ref{eqsmc1}).  All the numerical results are displayed for  $\mathbb{M}_{d}=\mathbb{S}_{2}=\left\lbrace\mathbf{x}\in\mathbb{R}^{3}: \ \|\mathbf{x}\|=1\right\rbrace,$  and, hence,  for $H=L^{2}(\mathbb{S}_{2},d\nu,\mathbb{R}),$ and   its complex version  $L^{2}(\mathbb{S}_{2},d\nu,\mathbb{C}).$

Note that the proposed methodology in Section \ref{secFANOVA} does not suppose a specific structure regarding the scalar time--varying   entries in the design matrix $\mathbf{X}.$  In our simulation study, we adopt the framework of one way ANOVA, assuming  the application of five treatments through time, respectively leading to five functional response spherical means, $\beta_{j}(\mathbf{x}),$ $\mathbf{x}\in\mathbb{S}_{2},$ $j=1,2,\dots,p.$ This example can be interpreted  in the field of change point analysis, since one can adopt  this framework to detecting changes  through time  in the spherical trend  of the observed functional time series (see, e.g.,  Dette \& Quanz \cite{DetteQuanz22}).

\begin{figure}[H]
\centering
\includegraphics[height=0.23\textheight, width=0.7\textwidth]{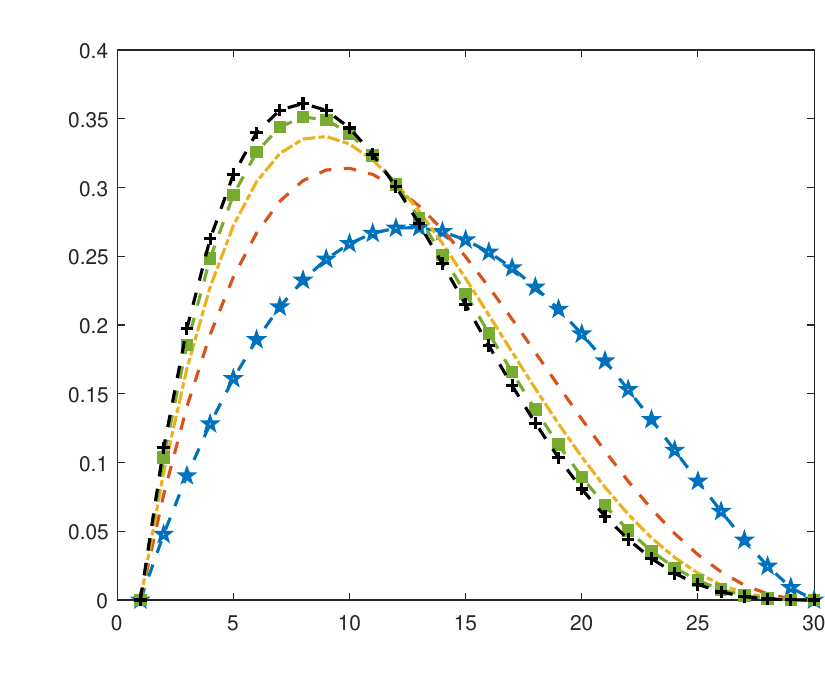}
\caption{Fourier coefficients  ${\beta}_{n,j},$ $n=1,2,\dots,30,$ for $j=1,2,3,4,5,$ are respectively displayed in blue star--dashed line,  red dashed line,  orange dotted line,  green squares-dashed line and  black cross--dashed line.}
\label{fcoefbetatruefanova}
\end{figure}

\begin{figure}[H]
\includegraphics[height=0.3\textheight, width=\textwidth]{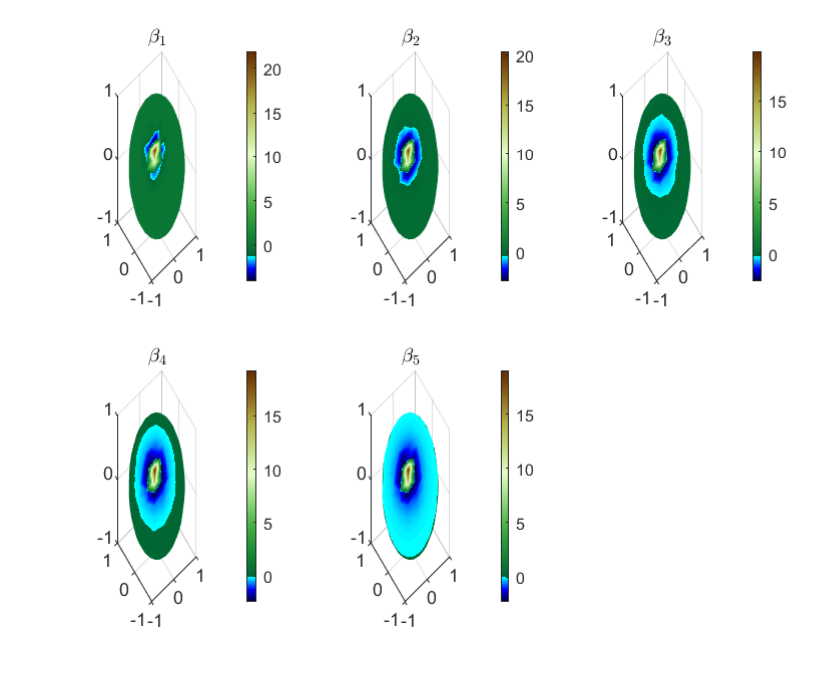}
\caption{The regression parameters $\beta_{j}(\mathbf{x}), \mathbf{x}\in \mathbb{S}_{2},j=1,2,\dots,5,$  projected into the direct sum of eigenspaces $\mathcal{H}_{n}$, $n=1,2,\dots,30,$ of the Laplace Beltrami operator on $L^{2}(\mathbb{S}_{2}, d\nu , \mathbb{R})$.}
\label{fsphbeta}
\end{figure}

 In the model  generations, truncation at discrete Legendre frequency $M=30$ has been achieved. The corresponding  Fourier coefficients of
the  $L^{2}(\mathbb{S}_{2},d\nu,\mathbb{R})$--valued parameters $\beta_{j},$  $j=1,2,\dots,p,$
 with respect the orthonormal basis of eigenfunctions of the Laplace Beltrami operator,
 are given by:

\begin{eqnarray}
\label{eceivbeta}
\beta_{n,j}&=&\frac{1}{6}\frac{x_{n}^{\alpha-1}(1-x_{n})^{\varrho_{j}-1}\Gamma(\alpha+\varrho_{j})}{\Gamma(\alpha)\Gamma(\varrho_{j})},\quad \alpha, \varrho_{j}>0\\
\alpha &=& 2,\quad  \varrho_{j}=\frac{5j}{j+1},\quad j=1,2,\dots,p,\quad  p=5\nonumber\\
x_{1}&=&0,\quad  x_{n}=x_{n-1}+\frac{1}{30-1},\quad  n=2,3,\dots,30.\nonumber
\end{eqnarray} \noindent  These Fourier coefficients are displayed in Figure \ref{fcoefbetatruefanova}. The projected regression parameters into the direct sum $\bigoplus_{n=1}^{30} \mathcal{H}_{n}$ of the eigenspaces $\mathcal{H}_{n},$ $n=1,\dots,30,$ of the Laplace Beltrami operator can also be seen in Figure \ref{fsphbeta}.

In this simulation study, we have considered the case where the error term $\varepsilon=\left\lbrace\varepsilon_{t},t\in\mathbb{Z}\right\rbrace$ obeys a multifractionally integrated SPHARMA($1,1$)
equation. Equivalently, the following state equation characterizes the behavior of $\varepsilon:$ For each  $t\in \mathbb{Z},$
\begin{eqnarray}
&&
(\mathcal{I}_{L^{2}(\mathbb{S}_{2},d\nu ,\mathbb{R})}-B)^{\mathcal{A}_{\theta_{0}}/2}(\boldsymbol{\Phi}_{p}(B)\varepsilon_{t})(\mathbf{x})=
\eta_{t}(\mathbf{x})
+
(\boldsymbol{\Psi}_{q}(B)\eta_{t})(\mathbf{x}),\quad \mathbf{x}\in \mathbb{S}_{2},\nonumber
\end{eqnarray}
\noindent
where $\mathcal{A}_{\theta_{0}}$ denotes the LRD integral operator with isotropic kernel $\mathcal{K}_{\mathcal{A}_{\theta_{0}}}$ satisfying $$\mathcal{K}_{\mathcal{A}_{\theta_{0}}}(\mathbf{x},\mathbf{y})=\sum_{n\in \mathbb{N}_{0}} \alpha (n,\theta_{0})
\sum_{j=1}^{\delta (n,d)}S_{n,j}^{d}(\mathbf{x})S_{n,j}^{d}(\mathbf{y}), \ \mathbf{x}, \mathbf{y}\in \mathbb{S}_{2},$$
\noindent with, as before,  $d=2,$ $\delta(n,d)=\delta(n,2)=2n+1.$ We have considered $l_{\alpha}(\theta_{0})\leq\alpha (n,\theta_{0})\leq L_{\alpha}(\theta_{0}),$  with  $l_{\alpha}(\theta_{0}), L_{\alpha}(\theta_{0})\in \left(0,1/2\right).$ $L^{2}(\mathbb{S}_{2},d\nu,\mathbb{R})$--valued strong white noise process $\eta=\left\lbrace\eta_{t}(\mathbf{x}), \ \mathbf{x}\in\mathbb{S}_{2}, \ t\in \mathbb{Z}\right\rbrace$ has  variance $\sigma^{2}_{\eta}=\sum_{n\in\mathbb{N}_{0}}\sigma^{2}_{n}.$  Furthermore,  $\boldsymbol{\Phi}_{p}(B)=1-\sum_{k=1}^{p}\Phi_{k}B^{k},$ and $\boldsymbol{\Psi}_{q}(B)=\sum_{l=1}^{q}\Psi_{l}B^{l}$  respectively denote  the autoregressive and moving average polynomial operators, with $\Phi_{k},k=1,2,\dots,p$ and $\Psi_{l},l=1,2,\dots,q$ being invariant positive self-adjoint bounded operators on $L^{2}(\mathbb{S}_{2},d\nu, \mathbb{R}).$ Each one of the above operators admits the diagonal series expansion:

\begin{eqnarray}
\Phi_{k}&=&\sum_{n\in \mathbb{N}_{0}}\lambda_{n}(\Phi_{k})\sum_{j=1}^{\delta (n,d)}S_{n,j}^{d}\otimes S_{n,j}^{d},\  k=1,\dots,p\nonumber\\
\Psi_{l}&=&\sum_{n\in \mathbb{N}_{0}}\lambda_{n}(\Psi_{l})\sum_{j=1}^{\delta (n,d)}S_{n,j}^{d}\otimes S_{n,j}^{d},\  l=1,\dots,q.\nonumber
\end{eqnarray}

\begin{figure}[H]
     \centering
     \begin{subfigure}[t]{0.49\textwidth}
         \centering
	\includegraphics[height=0.23\textheight, width=\textwidth]{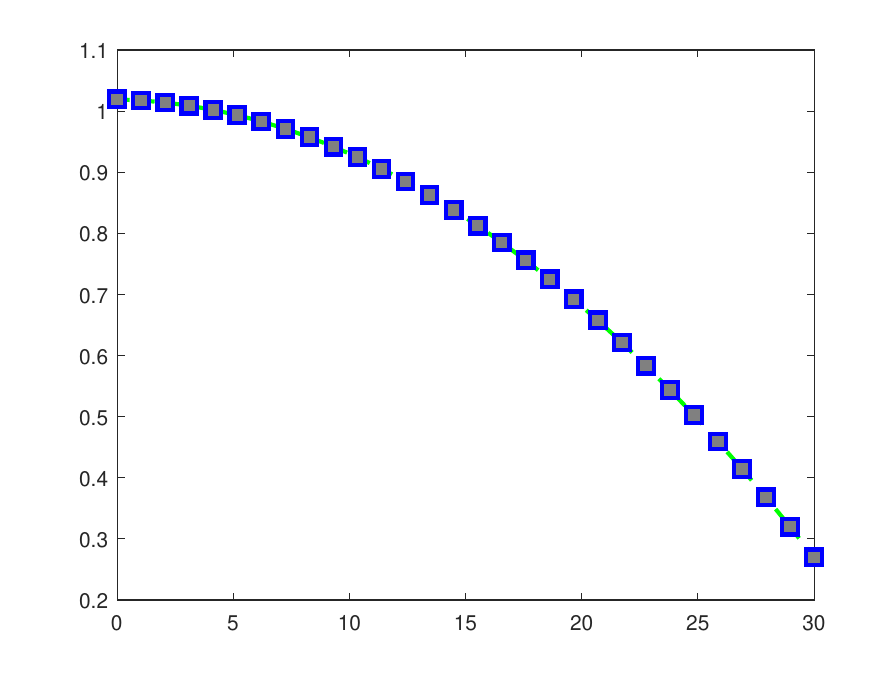}
\caption{}
     \end{subfigure}
     \hfill
     \begin{subfigure}[t]{0.49\textwidth}
         \centering
        \includegraphics[height=0.23\textheight, width=\textwidth]{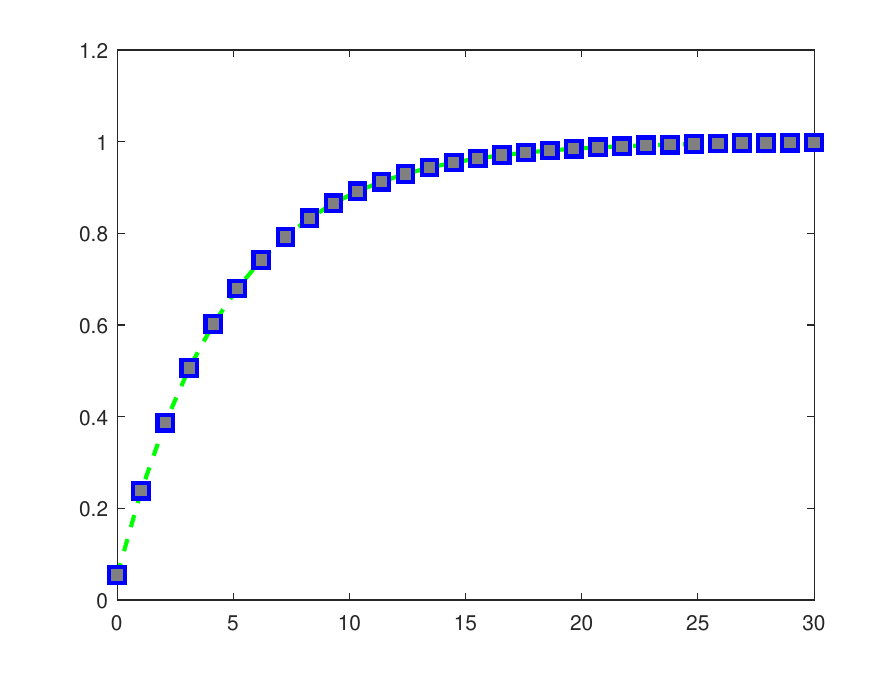}
         \caption{}
     \end{subfigure}
\caption{The first $30$ eigenvalues $\alpha (n,\theta_{0}),$   $n=1,\dots,30,$  of the  LRD operator $\mathcal{A}_{\theta_{0}}.$  The considered DPBS of eigenvalues is plotted at the left--hand side, and  the IPBS of eigenvalues  is plotted  at the right hand--side.}
\label{fsymLRD}
\end{figure}

The results displayed correspond to the generation of functional samples of size  $N=50,100,500$  of the multifractionally integrated  SPHARMA(1,1) process defining the error term in our regression model. In particular, the following parameter values have been considered in the generations of SPHARMA(1,1) process : $\sigma^{2}_{n}=(n+1)^{-3/2},$ $\lambda_{n}(\Phi_{k})=\left[0.7\left(n+\frac{1}{n}\right)\right]^{-3/2},$ $\lambda_{n}(\Psi_{l})=0.4\left(n+\frac{1}{n}\right)^{-3/2},$ $n=1,2,\dots,30,$ $k=l=1.$ Multifractionally integration has been achieved in the spectral domain in terms of the following parameterized
 eigenvalues $\{\alpha(n,\theta_{0}),\ n\in \mathbb{N}_{0}\}$ of the LRD operator $\mathcal{A}_{\theta_{0}}:$
\begin{eqnarray}
\label{eceivadec}
\alpha(n,\theta_{0})&=&\frac{\theta_{0,1}x_{n}^{2}+\theta_{0,2}x_{n}+\theta_{0,3}}{\theta_{0,4}}, \theta_{0}=(\theta_{0,1},\theta_{0,2},\theta_{0,3},\theta_{0,4}),\\
\theta_{0,1}&=&0.75, \theta_{0,2}=0.76, \theta_{0,3}=0.77,\nonumber \\
\theta_{0,4}&=&\sup_{i=1,2,\dots,100}f(x_{n},i),\nonumber\\
f(x_{n},i)&=&\frac{1}{100}\left(ix_{n}^{2}+(i+1)x_{n}+(i+2)\right),\nonumber\\
x_{1}&=&0,x_{n}=x_{n-1}+\frac{30}{29}, \ n=2,3,\dots,30,\nonumber\\
\alpha(n,\upsilon_{0})&=&1-\frac{1}{9\exp(\upsilon_{0,1}+\upsilon_{0,2}x_{n})}, \: \upsilon_{0}=(\upsilon_{0,1},\upsilon_{0,2})=(1,1),\label{eceivainc}\\
x_{1}&=&-\pi,x_{n}=x_{n-1}+\frac{2\pi}{30-1},\: n=2,3,\dots,30,\nonumber
\end{eqnarray}
\noindent where  we have considered a Decreasing Positive Bounded Sequence (DPBS) and an Increasing Positive Bounded Sequence (IPBS) of eigenvalues  $\left\{\alpha(n,\theta_{0}), \ n\in \mathbb{N}_{0}\right\}$ and   $\left\{\alpha(n,\upsilon_{0}),\ n\in \mathbb{N}_{0}\right\}$ of the
LRD operator, respectively defined  in equations (\ref{eceivadec}) and (\ref{eceivainc})  (see also Figure \ref{fsymLRD}).

 The values  $l_{\alpha}=1.0192$ and $L_{\alpha }=0.2699$ have been considered in  DPBS, and $l_{\alpha }=0.0541$ and $L_{\alpha }=0.9982$ in IPBS of eigenvalues  of the LRD  operator.
 Note that the dominant eigenspace of the Laplace Beltrami operator where the largest dependence range in time is displayed by the projected error term  under the DPBS corresponds to the  first eigenspace plotted, while, in the IPBS case, the projected process displays the strongest LRD in the last eigenspace  of the Laplace Beltrami operator plotted in Figure \ref{fsymLRD}.

\subsection{Spherical multiscale residual  analysis under  totally specified model}
\begin{figure}[H]
\includegraphics[width=\textwidth, height=0.35\textheight]{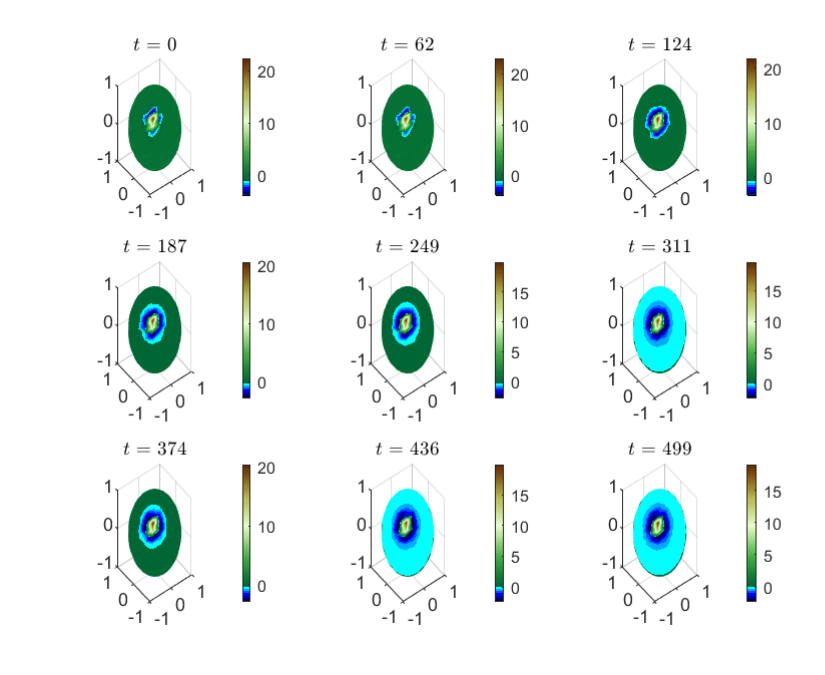}
\caption{REM $\hat{E}[Y_{t}(\mathbf{x})],$ $\mathbf{x}\in \mathbb{S}_{2},$  based on $R=100$ repetitions,   at times   $t=0,62,124,187,249,311,374,436,499,$ projected into the direct sum of the eigenspaces $\mathcal{H}_{n}$, $n= 1,2,\dots,30$ of the Laplace Beltrami operator on $L^{2}(\mathbb{S}_{2}, d\nu, \mathbb{R})$, under    DPBS  of eigenvalues  of LRD  operator. The functional sample size generated is $T=500$.}
\label{fSPHARMA11YMeandecre}
\end{figure}

\begin{figure}[H]
\includegraphics[width=\textwidth, height=0.35\textheight]{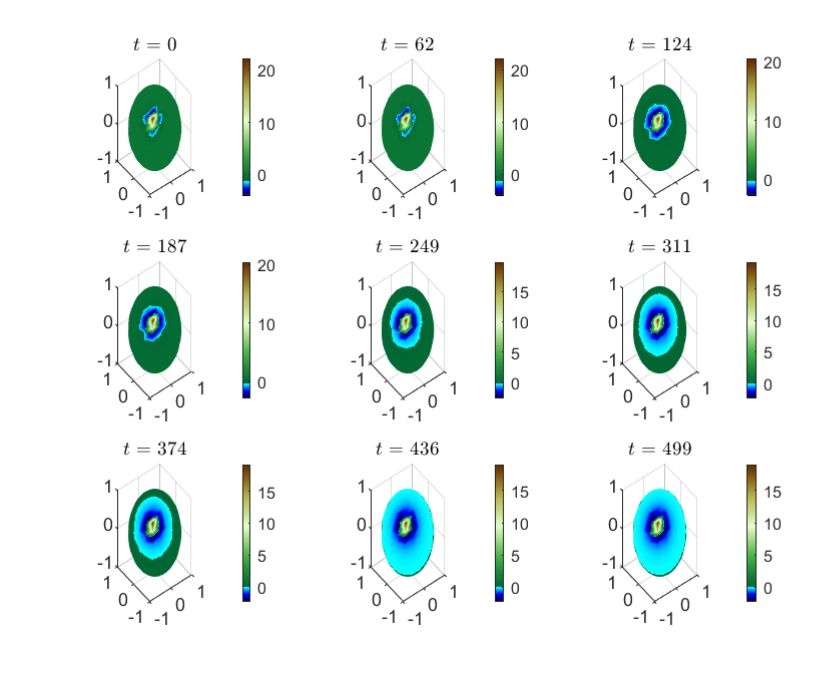}
\caption{RTPEM $\hat{E}\left[\hat{Y}_{t}(\mathbf{x})\right],$  $\mathbf{x}\in \mathbb{S}_{2},$   based on $R=100$ repetitions, at times $t=0,62,124,187,249,311,374,$ $436,499,$    projected into the direct sum of the eigenspaces $\mathcal{H}_{n}$, $n= 1,2,\dots,30$ of the Laplace Beltrami operator on $L^{2}(\mathbb{S}_{2}, d\nu, \mathbb{R})$, under    DPBS  of eigenvalues  of LRD  operator. The functional sample size generated is $T=500$.}
\label{fSPHARMA11YestMeandecre}
\end{figure}

The generations displayed in Figures  \ref{fSPHARMA11YMeandecre}  and \ref{fSPHARMA11YestMeandecre} correspond to the case of a  DPBS of eigenvalues of the LRD operator. Specifically,  the  empirical mean of the response  (REM), based on $R=100$ repetitions of a functional sample of size $N=500,$     projected into the direct sum $\bigoplus_{n=1}^{30}\mathcal{H}_{n}$ of the eigenspaces $\mathcal{H}_{n},$ $n=1,\dots,30,$ of the Laplace Beltrami operator  is showed in Figure \ref{fSPHARMA11YMeandecre} at times $t=0,62,124,187,249,311,374,436,499.$  The empirical mean  (RTPEM) (based on the same $R=100$ repetitions) of the theoretical predictor, projected into the direct sum $\bigoplus_{n=1}^{30}\mathcal{H}_{n},$
is shown in Figure \ref{fSPHARMA11YestMeandecre} at times $t=0,62,124,187,249,311,374,436,499.$

\begin{figure}[H]
\centering
\includegraphics[width=\textwidth, height=0.35\textheight]{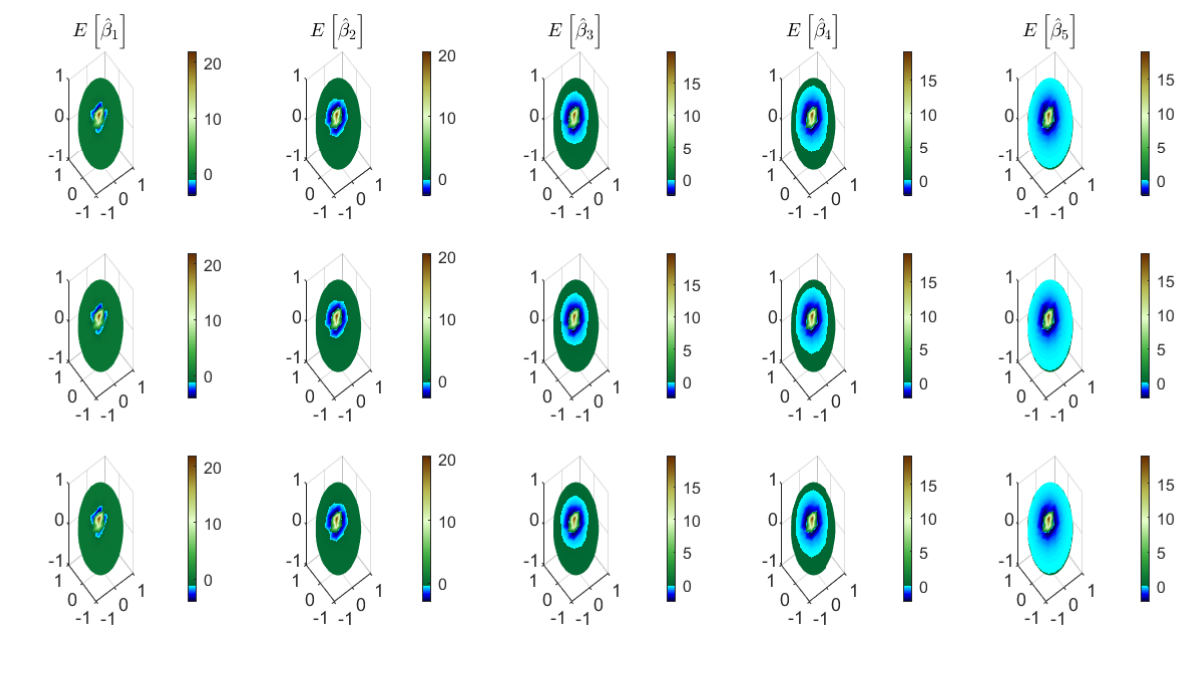}
\caption{Under DPBS of eigenvalues of the LRD operator, the projected empirical mean $\hat{E}\left[\hat{\beta}_{j}(\mathbf{x})\right],$   $\mathbf{x}\in\mathbb{S}_{2},$ of the computed values of the  GLS parameter estimator, based on $R=100$ repetitions,  are showed for $j=1,2,3,4,5.$   The results for functional sample sizes $N=50,100,500$ are  respectively presented by rows.}
\label{fSPHARMA11meanbetaestdecre}
\end{figure}

\begin{figure}[H]
\centering
\includegraphics[width=\textwidth, height=0.35\textheight]{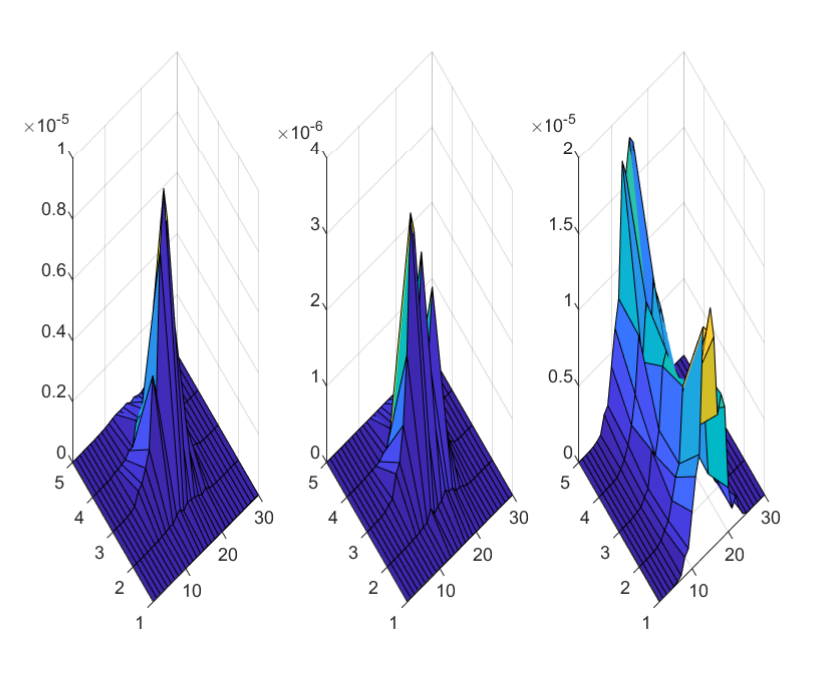}
\caption{Under  DPBS of eigenvalues of the LRD operator, the Empirical Mean Quadratic Errors (EMQEs) $\widehat{E}\left[(\widehat{\beta}_{n,j}-\beta_{n,j})^2\right], n= 1,2,\dots,30,$ $j=1,2,\dots,5,$ based on $R=100$ repetitions, for functional sample sizes $N=50,100,500$ are respectively plotted from the left--hand--side to the right--hand--side.}
\label{fARMA11MECBETAAdecre}
\end{figure}

\begin{figure}[H]
\centering
\includegraphics[width=\textwidth, height=0.3\textheight]{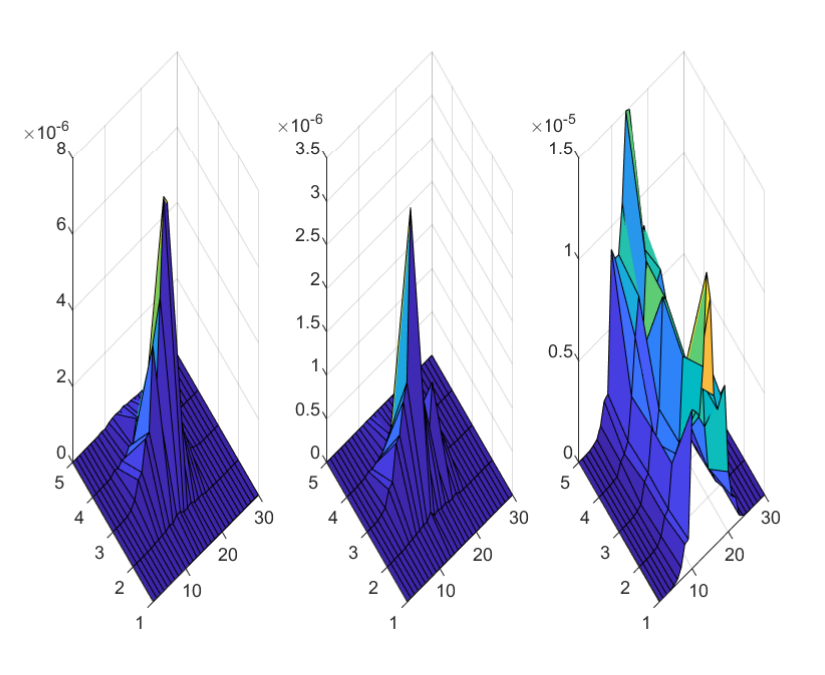}
\caption{Under  IPBS of eigenvalues of the LRD operator, the  EMQEs  $\widehat{E}\left[(\widehat{\beta}_{n,j}-\beta_{n,j})^2\right],  n= 1,2,\dots,30,$ $j=1,2,\dots,5,$   based on $R=100$ repetitions, for functional sample sizes $N=50,100,500$ are respectively plotted from the left--hand--side to the right--hand--side.}
\label{fARMA11MECBETAAcre}
\end{figure}

Figure \ref{fSPHARMA11meanbetaestdecre} illustrates unbiasedness of $\widehat{\beta}_{j}(\mathbf{x}),$ $\mathbf{x}\in\mathbb{S}_{2},$  $j=1,2,3,4,5,$ whose empirical  means based on $R=100$  repetitions under    DPBS  of eigenvalues  of LRD  operator are plotted in this figure for  functional sample sizes $N=50,100,500.$
The slow temporal decay of time--varying eigenvalue $B_{n}(t)$ as $t$ increases acts as a regularizer, regarding the singular behavior of the linear filter defining  the GLS parameter estimator in   (\ref{ecbetamcg})  as $n$ increases  for a fixed time $t.$ Note that this singular behavior is induced
by  $\Lambda_{n}^{-1}$ as  $n$ increases due to the fast decay of the pure point spectra of the elements of the trace covariance operator family.
This regularization effect, which becomes stronger when the functional sample size increases,   can be observed at least up to discrete Legendre frequency $n=M=30,$ being  more pronounced  at  low  and high discrete Legendre frequencies under DPBS and  IPBS  of eigenvalues  of the  LRD  operator, respectively. Thus, the highest  empirical mean quadratic errors  are observed at  spherical  scales corresponding to  intermediate Legendre frequencies  $n\in [10,20]$  (see, Figures \ref{fARMA11MECBETAAdecre} and \ref{fARMA11MECBETAAcre}).

\begin{figure}[H]
\centering
\includegraphics[width=\textwidth, height=0.3\textheight]{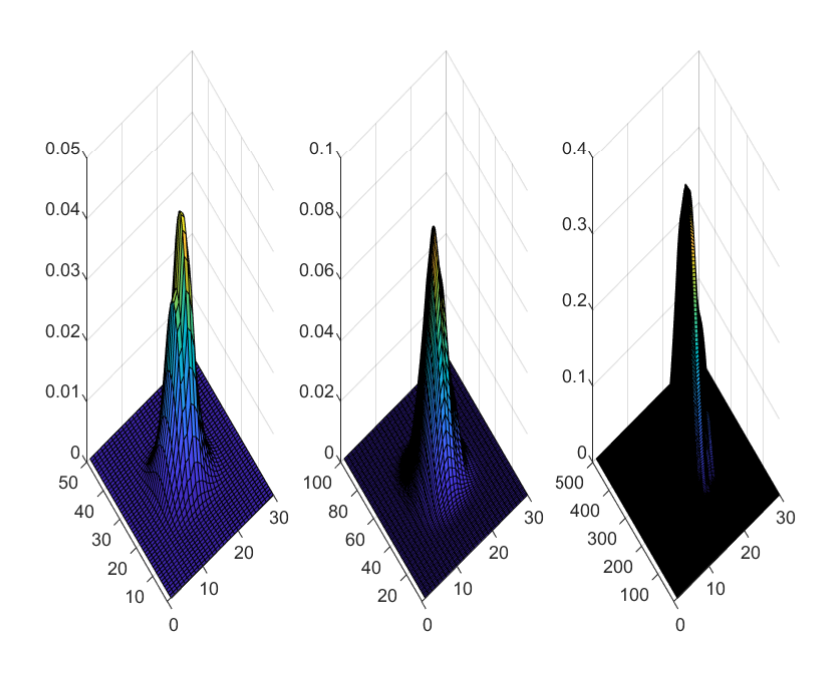}
\includegraphics[width=\textwidth, height=0.3\textheight]{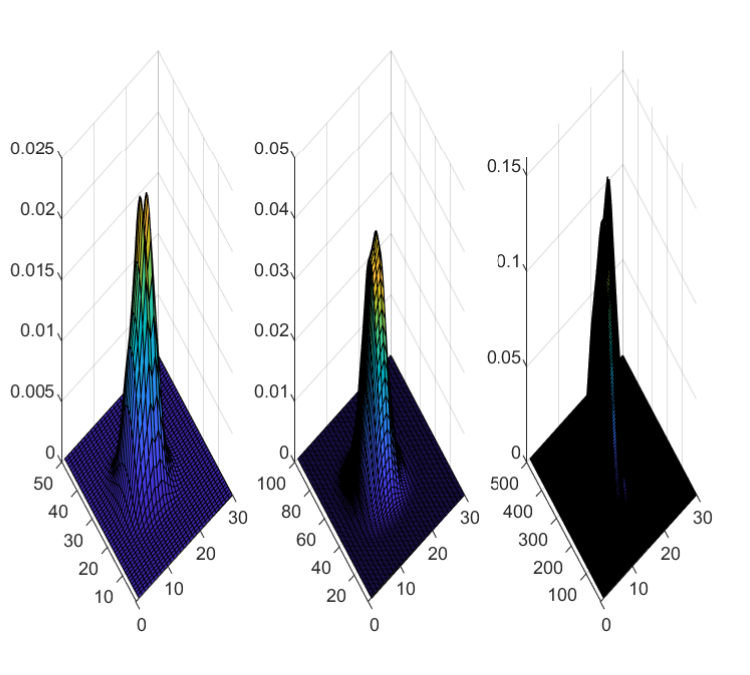}
\caption{Projected EMQEs of the theoretical predictor projections $\widehat{E}\left[(\widehat{Y}_{\mathbf{n}}(t)-Y_{\mathbf{n}}(t))^2\right],$    into $30$ eigenspaces of the Laplace Beltrami operator, for  $t=1,2,\dots,N$,
based on $R=100$ repetitions,
for functional sample sizes $N=50,100,500$ from the left--hand--side to the right--hand--side, respectively, under the IPBS (at the top), and under DPBS (at the bottom) of eigenvalues of the LRD operator.}
\label{fEMQEARMA11YTPcrecc}
\end{figure}

\begin{figure}[H]
\centering
\includegraphics[width=\textwidth, height=0.35\textheight]{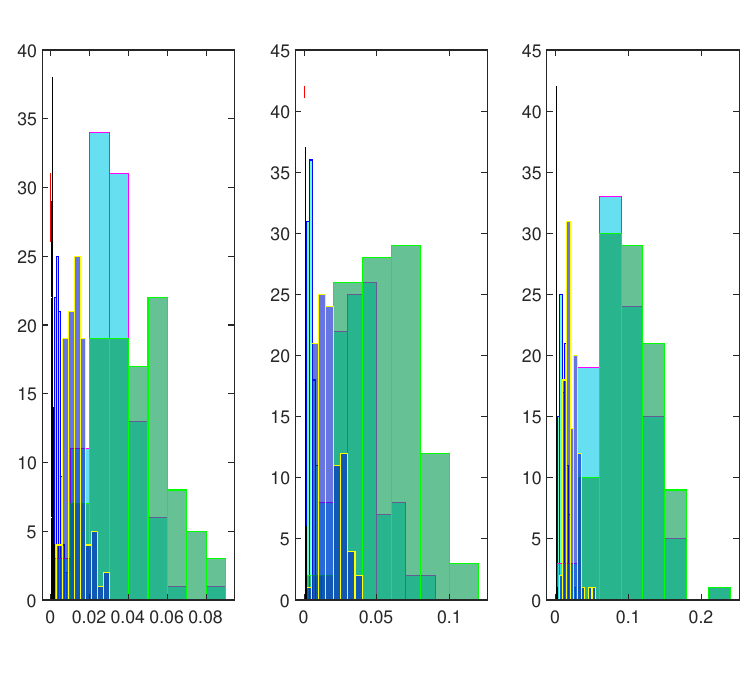}
\caption{Empirical distribution, based on $R=100$ repetitions, of $||\hat{Y}_{\mathbf{n}}(t)-Y_{\mathbf{n}}(t)||_{L^{1}([1,N])},$ for projections into Laplace--Beltrami eigenspaces $\mathcal{H}_{n},$ $n=1,5,10,15,20,25,30,$ for  functional sample sizes $T=50,100,500$ from the left--hand side to the right--hand side, respectively, under the DPBS of eigenvalues of the LRD operator.}
\label{fAETPARMA11decre}
\end{figure}

\begin{figure}[H]
\centering
\includegraphics[width=\textwidth, height=0.35\textheight]{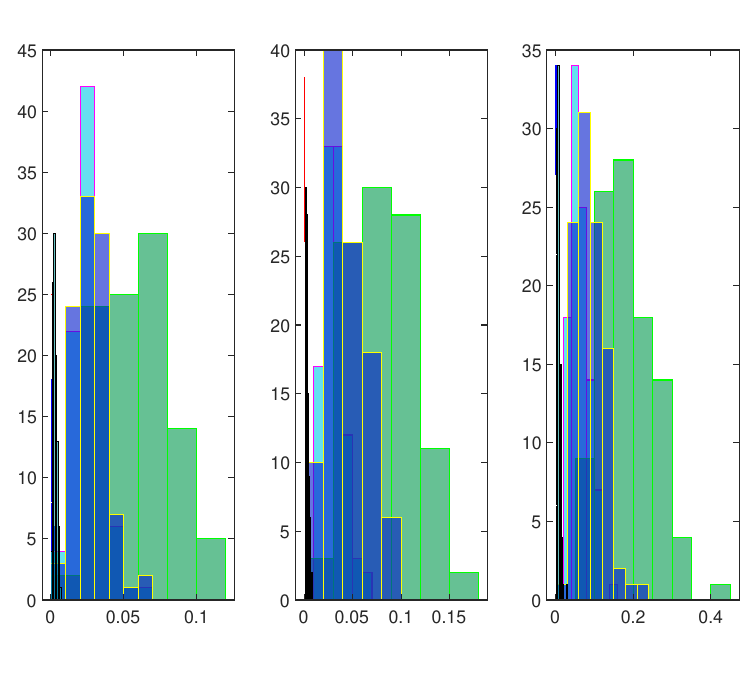}
\caption{Empirical distribution, based on $R=100$ repetitions, of $||\hat{Y}_{\mathbf{n}}(t)-Y_{\mathbf{n}}(t)||_{L^{1}([1,N])},$ for projections into Laplace--Beltrami eigenspaces $\mathcal{H}_{n},$ $n=1,5,10,15,20,25,30,$ for  functional sample sizes $T=50,100,500$ from the left--hand side to the right--hand side, respectively, under the IPBS of eigenvalues of the LRD operator.}
\label{fAETPARMA11cre}
\end{figure}

The incorporation of the time--varying coefficients defined by the entries of the design matrix leads to the definition of the theoretical predictor in equation (\ref{thpred}). The  temporal  pointwise values  of the functional Empirical Mean Quadratic Errors (EMQEs), based on $R=100$ repetitions,   associated with the  theoretical response predictor reflect the effect of the spherical--scale--varying strong--dependence  in time displayed by the response. Hence, the same conclusions   follow as in the  GLS parameter estimator (see Figure \ref{fEMQEARMA11YTPcrecc}). Note that, in the     analysis performed in terms of the empirical distribution of the $L^{1}$ norms of the  response prediction errors at each eigenspace $\mathcal{H}_{n}$  of the Laplace Beltrami operator, for $n=1,\dots,30,$ a slighter LRD effect is observed  in all spherical scales analyzed under DPBS than under IPBS of  eigenvalues of the LRD operator, as one can observe from the supports and modes of the empirical distributions plotted for the  functional sample sizes $N=50, 100, 500$ in Figures   \ref{fAETPARMA11decre} and \ref{fAETPARMA11cre}.

\subsection{Results under misspecified model}
\begin{figure}[H]
\centering
\includegraphics[width=\textwidth, height=0.45\textheight]{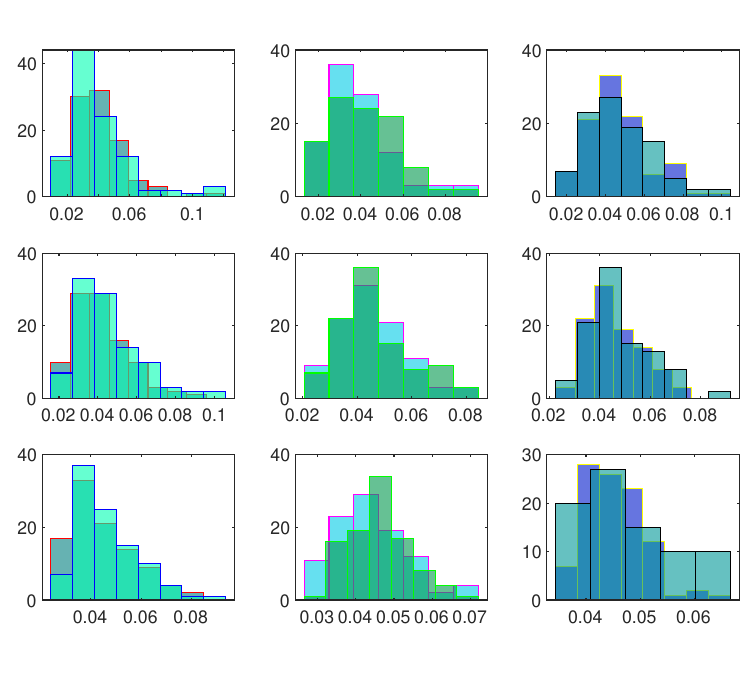}
\caption{Under DPBS of eigenvalues of the LRD operator, the empirical distribution, based on $R=100$ repetitions, of $||f_{n,\widehat{\theta}_{N}}(\cdot )-f_{n,\theta_{0}}(\cdot )||_{L^{1}([-\pi,\pi])}$ is showed. The results are displayed for
projections into   $\mathcal{H}_{n},$ for $n= 5-10$ at the left hand--side, for $n=15-20$ at the center, and for $n=25-30$ at the right--hand side.
By rows, the  empirical distributions plotted correspond to the functional sample sizes $N=50$ (first row), $N=100$ (second row) and $N=500$ (third row). }
\label{fARMA11histMCEdecre}
\end{figure}

\begin{figure}[H]
\centering
\includegraphics[width=\textwidth, height=0.45\textheight]{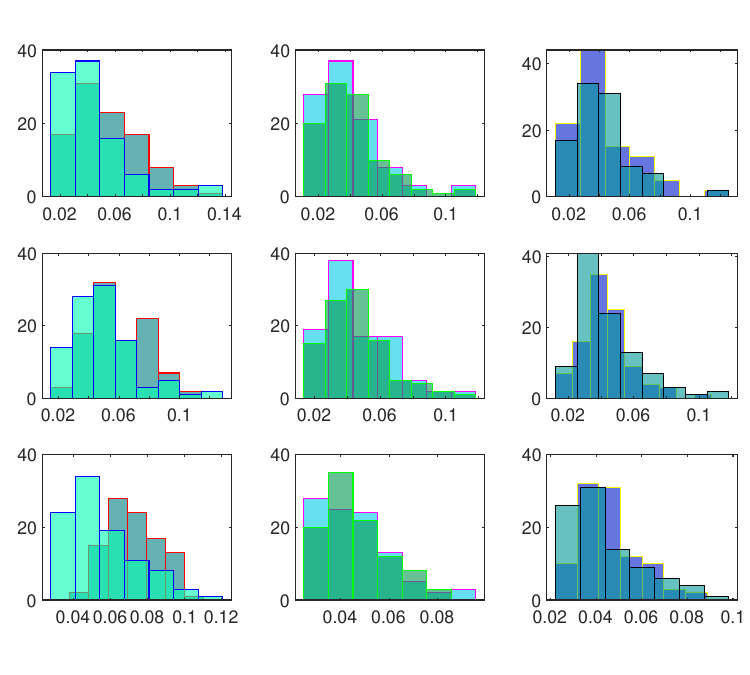}
\caption{Under IPBS  of eigenvalues of the LRD operator,
the empirical distribution, based on $R=100$ repetitions, of $||f_{n,\hat{\theta}_{N}}(\cdot )-f_{n,\theta_{0}}(\cdot )||_{L^{1}([-\pi,\pi])}$ is showed. The results are displayed for
projections into   $\mathcal{H}_{n},$ for $n= 5-10$ at the left hand--side, for $n=15-20$ at the center, and for $n=25-30$ at the right--hand side.
By rows, the  empirical distributions plotted correspond to the functional sample sizes $N=50$ (first row), $N=100$ (second row) and $N=500$ (third row).}
\label{fARMA11histMCEcre}
\end{figure}

The numerical results displayed in this section are obtained from the plug--in estimation of $\boldsymbol{\beta},$ based on the minimum contrast estimation of the parameterized frequency--varying eigenvalues of the elements of the spectral density operator family, as given in Section
\ref{secmiss}. Figures \ref{fARMA11histMCEdecre} and \ref{fARMA11histMCEcre} display the empirical distribution of the $L^{1}$--norm, based on $R=100$ repetitions, of the empirical absolute functional errors, associated with minimum contrast estimation from  the generated functional samples of sizes $N=50,100, 500.$

\begin{figure}[H]
\includegraphics[width=\textwidth, height=0.35\textheight]{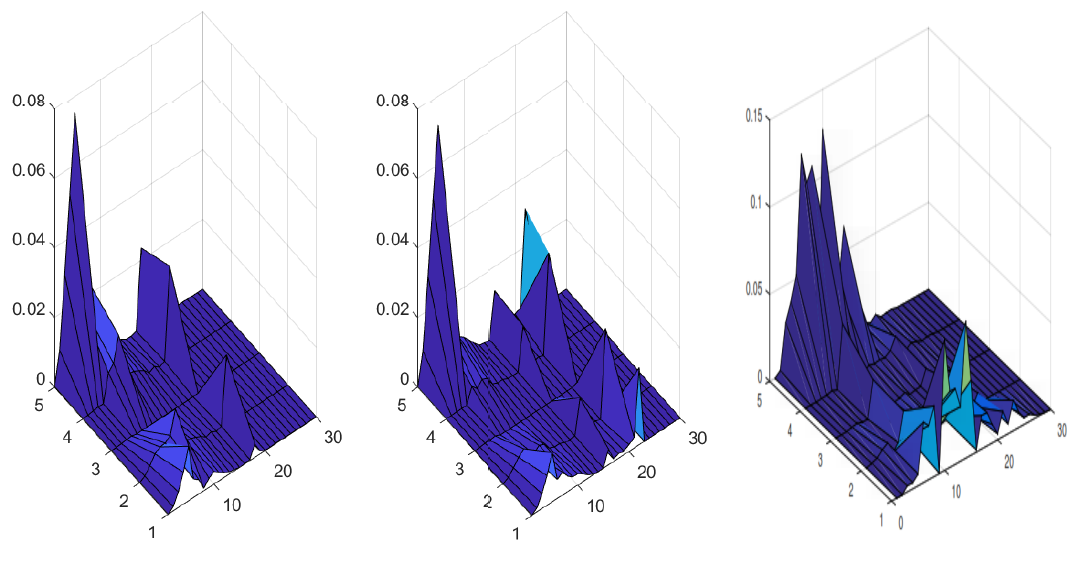}
\caption{Under DPBS of  eigenvalues of the LRD operator, the EMQE $\widehat{E}\left[(\widehat{\beta}_{n,j}-\beta_{n,j})^2\right], n= 1,2,\dots,30,$ $j=1,2,\dots,5$,  based on $R=100$ repetitions, are displayed for functional sample sizes $N=50$  at the left hand--side, $N=100$ at center,  and $N=500$    at the right hand--side.}
\label{fEMQEARMA11betaestmissdecre}
\end{figure}

\begin{figure}[H]
\includegraphics[width=\textwidth, height=0.35\textheight]{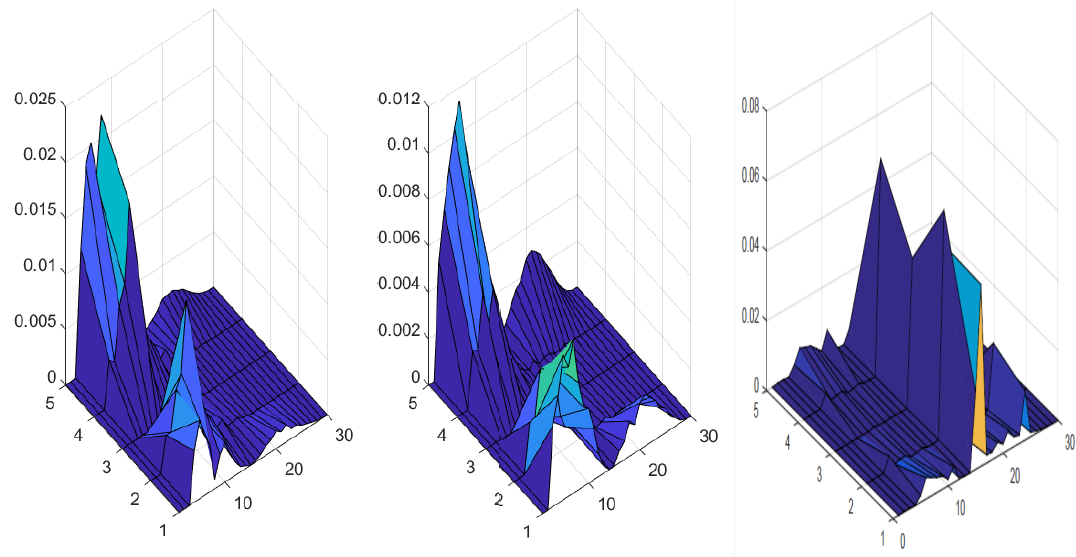}
\caption{Under IPBS of eigenvalues  of the LRD operator, the  EMQE $\widehat{E}\left[(\widehat{\beta}_{n,j}-\beta_{n,j})^2\right],$ $n= 1,2,\dots,30,$ $j=1,2,\dots,5,$ based on $R=100$ repetitions, are displayed for functional sample sizes $N=50$  at the left hand--side, $N=100$ at center,  and $N=500$    at the right hand--side.}
\label{fEMQEARMA11betaestmisscre}
\end{figure}

\begin{figure}[H]
\centering
\includegraphics[width=\textwidth, height=0.35\textheight]{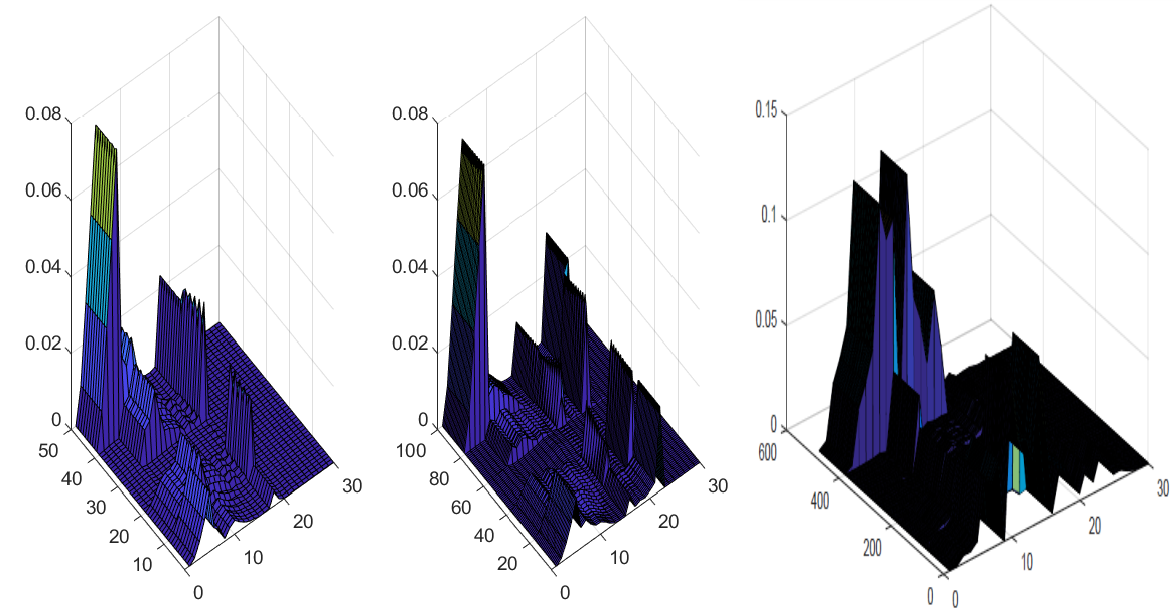}
\caption{Under DPBS of  eigenvalues of the LRD operator,   the EMQE $\widehat{E}\left[(\widehat{Y}_{\mathbf{n}}(t)-Y_{\mathbf{n}}(t))^2\right],$ $n= 1,2,\dots,30,$  $t=1,2,\dots,N,$ based on $R=100$ repetitions, are displayed for functional sample sizes $N=50$  at the left hand--side, $N=100$ at center,  and $N=500$    at the right hand--side.}
\label{fEMQEARMA11YTPdecre}
\end{figure}

\begin{figure}[H]
\centering
\includegraphics[width=\textwidth, height=0.35\textheight]{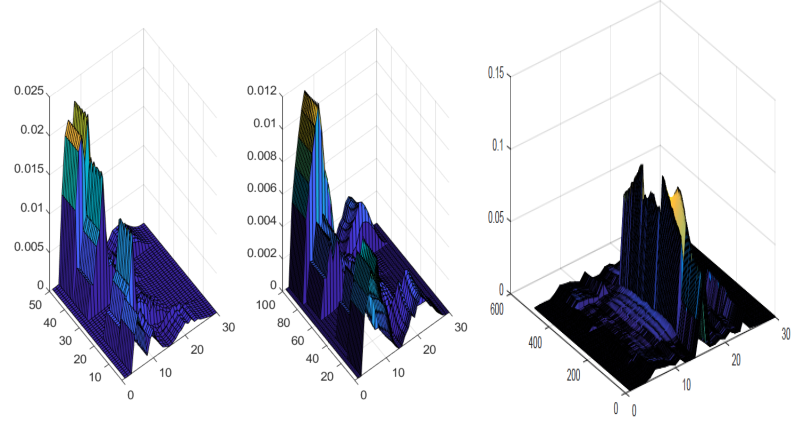}
\caption{Under IPBS  of eigenvalues of the LRD operator, the EMQE $\widehat{E}\left[(\widehat{Y}_{\mathbf{n}}(t)-Y_{\mathbf{n}}(t))^2\right],$ $n= 1,2,\dots,30,$ $t=1,2,\dots,N,$ based on $R=100$ repetitions, are displayed for functional sample sizes $N=50$  at the left hand--side, $N=100$ at center,  and $N=500$    at the right hand--side.}
\label{fEMQEARMA11YTPcre}
\end{figure}

 The resulting spherical scale dependent EMQEs, associated  with the plug--in estimation of the Fourier coefficients of the  functional regression parameter vector   $\boldsymbol{\beta }$ from its projection into the  eigenspaces $\mathcal{H}_{n},$ $n=1,\dots,30,$ of the Laplace Beltrami operator are then computed (see Figures  \ref{fEMQEARMA11betaestmissdecre} and \ref{fEMQEARMA11betaestmisscre}).  The corresponding temporal pointwise values of the response prediction EMQEs are also displayed by spherical scales (corresponding to projection into $\mathcal{H}_{n},$ $n=1,\dots,30,$) for functional sample sizes $N=50,100, 500,$  based on $R=100$ repetitions  (see Figures \ref{fEMQEARMA11YTPdecre}  and \ref{fEMQEARMA11YTPcre}). All the above numerical results are visualized under  DPBS and  IPBS  of eigenvalues of  the LRD operator.

Regarding  the empirical distribution of $||f_{n,\widehat{\theta}_{N}}(\cdot )-f_{n,\theta_{0}}(\cdot )||_{L^{1}([-\pi,\pi])},$ plotted in Figures  \ref{fARMA11histMCEdecre}  and \ref{fARMA11histMCEcre}, under DPBS and IPBS of eigenvalues of the LRD operator, respectively,
one  can  observe that   IPBS  of  eigenvalues of the LRD operator slightly enlarges  the supports of the empirical distributions with respect to DPBS of  eigenvalues of the LRD operator. Also, more asymmetric patterns are observed, in some spherical  scales, and, in some cases,  closer to  bimodality (see, e.g., right--hand side plot at the  third row of Figure \ref{fARMA11histMCEcre}). Thus, under DPBS of  eigenvalues of the LRD operator, more symmetry is observed in the empirical distribution patterns     at most of the spherical scales, and  the empirical distribution  supports reduce faster, when the functional sample size increases. These differences are stronger at the highest resolution levels displayed (i.e., when para meter $n$ increases).

The plotted results in Figures  \ref{fARMA11histMCEdecre}  and \ref{fARMA11histMCEcre} on the effect of parameters $n$ (spherical scale),
 and  $N$ (functional sample size) on the empirical distribution of
 the statistics $||f_{n,\widehat{\theta}_{N}}(\cdot )-f_{n,\theta_{0}}(\cdot )||_{L^{1}([-\pi,\pi])}$ lead to a stronger level of singularity in the
 estimated temporal  functional spectrum under DPBS   than under IPBS  of  eigenvalues of  the LRD operator  at coarser spherical scales ($n$ small), due to a stronger smoothing effect of the latest (the IPBS of  eigenvalues of  the LRD operator  scenario)   at the zero frequency pole, for such coarser spherical scales. While at high resolution levels ($n$ large), the opposite effect is observed regarding the stronger smoothing at high resolution levels  of zero frequency pole under
DPBS of eigenvalues  of the LRD operator. These differences between both analyzed LRD operator eigenvalues scenarios  are more pronounced when the functional sample size $N$  increases, as it can be observed in Figures \ref{fEMQEARMA11betaestmissdecre}--\ref{fEMQEARMA11YTPcre}, regarding the highest  EMQEs values displayed at coarser and high resolution levels in the sphere, associated with the plug--in estimation of $\boldsymbol{\beta },$ and the corresponding prediction of the response $\mathbf{Y}.$
\section{Conclusions and open research lines}
\label{secconclusions}

This paper opens a new research line within the context of  multiple functional regression analysis from strong--correlated in time functional data in a manifold $\mathbb{M}_{d}.$ Particularly, the framework of   connected  and compact two--point homogeneous spaces is adopted. The  formulated multiple functional regression model, with functional response, functional regression parameters and time--dependent scalar covariates, goes beyond the  assumptions of weak--dependent, and the Euclidean setting usually adopted in the current literature in functional regression. The asymptotic properties of the
GLS estimator of the  $\mathbb{M}_{d}$--supported functional regression parameter vector $\boldsymbol{\beta },$ when the second--order structure of the functional error term is known and unknown,    are established. Indeed, these properties are obtained by applying the results derived  in Ruiz--Medina, Miranda \& Espejo \cite{RuizMD18}.
In the GLS estimation of the  $\mathbb{M}_{d}$--supported functional regression parameter vector $\boldsymbol{\beta },$ the simulation study undertaken shows that, under totally specified model, i.e., when the second--order structure of the functional error term $\varepsilon $ is known, the velocity decay and distribution of the singular values of the elements of the covariance operator family play a crucial role, in relation to the regularization effect induced by the slow decay in time of the covariance function. Note that this regularization effect also depends on the distribution  of the eigenvalues of the LRD operator (in our analysis it depends  on the DPBS and  IPBS of eigenvalues of the LRD operator considered).
The regularization effect acts
 on the inverse matrix sequence
$\left\{\Lambda_{n}^{-1},\ n\in \mathbb{N}_{0}\right\},$ involved in the linear filter defining the GLS parameter estimator.
 When the functional spectral--based plug--in estimation of $\boldsymbol{\beta }$  is achieved under misspecified model, the increasing of the functional sample size leads to a more dense sampling at a neighborhood of the zero frequency pole, increasing the magnitude of the errors of the computed minimum contrast parameter estimator of the LRD operator. This effect is mitigated in part, in our study, by the stronger  smoothing effect at spherical coarser scales under IPBS of eigenvalues of the
  LRD operator, while the opposite effect is observed under  DPBS of eigenvalues of the
  LRD operator at such scales, with a stronger smoothing effect at high resolution levels in the sphere.

In a subsequent paper we address significance testing under this strong--dependence functional data scenario in the analyzed manifold  multiple functional regression model (see, e.g., Garc\'{\i}a--Portugu\'es,  Gonz\'alez--Manteiga \& Febrero--Bande \cite{Febrero}, when i.i.d functional data are analyzed, and Ruiz--Medina \cite{RuizMedina16} and
\'{A}lvarez--Li\'{e}bana \&   Ruiz--Medina \cite{Alvarez17}  for  weak--dependent functional data in the time series context).

\section*{Acknowledgements}\label{acknowledgements}

This work has been supported in part by projects
MCIN/ AEI/\linebreak PID2022-142900NB-I00,  and  CEX2020-001105-M MCIN/ AEI/10.13039/501100011033.
 
 \end{document}